\documentclass[graybox]{svmult}

\usepackage{type1cm}        %
\usepackage{makeidx}         %
\usepackage{graphicx}        %
\usepackage{multicol}        %
\usepackage[bottom]{footmisc}%

\usepackage{newtxtext}       %
\usepackage[varvw]{newtxmath}       %
\usepackage{microtype}
\usepackage{nicefrac}
\usepackage{booktabs}
\usepackage{subcaption}

\usepackage{tikz}
\usepackage{tikz-network}
\usetikzlibrary{patterns}
\usepackage{pgfplots}
\pgfplotsset{compat=1.18}
\definecolor{nvidiaGreen}{RGB}{118,185,0}
\usetikzlibrary{math} %
\usepackage[colorlinks]{hyperref}
\usepackage{listings}
\usepackage{comment}
\lstdefinestyle{fortran}{%
  backgroundcolor=\color{yellow!10},%
  basicstyle=\small\ttfamily,%
  identifierstyle=\color{black},%
  keywordstyle=\color{blue},%
  keywordstyle={[2]\color{cyan}},%
  keywordstyle={[3]\color{olive}},%
  stringstyle=\color{teal},%
  commentstyle=\itshape\color{orange},%
}
\usepackage{siunitx}

\DeclareMathOperator*{\argmax}{arg\,max}

\usepackage{todonotes}

\makeindex             %

\begin{document}

\title*{Parallel matching-based AMG preconditioners for elliptic equations discretized by IgA}
\titlerunning{Parallel AMG preconditioners for IgA} %
\author{Pasqua D'Ambra\orcidID{0000-0003-2047-4986} and\\ Fabio Durastante\orcidID{0000-0002-1412-8289} and\\ Salvatore Filippone\orcidID{0000-0002-5859-7538}}
\authorrunning{P. D'Ambra, F. Durastante, and S. Filippone}
\institute{Pasqua D'Ambra \at National Research Council of Italy (CNR), Institute for Applied Computing ``Mauro Picone'' (IAC), \email{pasqua.dambra@cnr.it}
\and Fabio Durastante \at University of Pisa, Largo Bruno Pontecorvo, 5 -- 56127 Pisa and IAC-CNR \email{fabio.durastante@unipi.it}, \and
Salvatore Filippone \at University of Rome ``Tor Vergata'', Via Politecnico 1 -- 00133 Roma and IAC-CNR \email{salvatore.filippone@uniroma2.it}
}
\maketitle

\abstract*{Isogeometric analysis (IgA) offers enhanced approximation capabilities for the discretization of elliptic boundary-value problems, yet it results in large, sparse, and increasingly ill-conditioned linear systems due to higher interconnectivity among degrees of freedom. In particular, the discretization with tensor-product B-splines or NURBS of degree $p$ on a mesh with $n$ elements per parametric direction leads to symmetric positive-definite systems of the form $K\mathbf{u} = \mathbf{F}$, where the matrix bandwidth and condition number scale unfavorably with both $p$ and spatial dimension $d$. To address the computational challenges posed by such systems, especially in three-dimensional or high-order scenarios, Krylov subspace methods with specialized preconditioners become essential. This paper investigates the efficacy of algebraic multigrid (AMG) preconditioners tailored for IgA-based discretizations, with a focus on performance in modern high-performance computing (HPC) environments. Leveraging the Parallel Sparse Computation Toolkit (\texttt{PSCToolkit}), we explore distributed-memory and GPU-accelerated strategies for solving large-scale problems. The study assesses algorithmic efficiency and scalability across a range of benchmark tests. The results demonstrate that AMG preconditioners can achieve robust and scalable performance, confirming their potential as practical solvers for large IgA systems in engineering and scientific applications.}

\abstract{Isogeometric analysis (IgA) offers enhanced approximation capabilities for the discretization of elliptic boundary-value problems, yet it results in large, sparse, and increasingly ill-conditioned linear systems due to higher interconnectivity among degrees of freedom. In particular, the discretization with tensor-product B-splines or NURBS of degree $p$ on a mesh with $n$ elements per parametric direction leads to symmetric positive-definite systems of the form $K\mathbf{u} = \mathbf{F}$, where the matrix bandwidth and condition number scale unfavorably with both $p$ and spatial dimension $d$. To address the computational challenges posed by such systems, especially in three-dimensional or high-order scenarios, Krylov subspace methods with specialized preconditioners become essential. This paper investigates the efficacy of algebraic multigrid (AMG) preconditioners tailored for IgA-based discretizations, with a focus on performance in modern high-performance computing (HPC) environments. Leveraging the Parallel Sparse Computation Toolkit (\texttt{PSCToolkit}), we explore distributed-memory and GPU-accelerated strategies for solving large-scale problems. The study assesses algorithmic efficiency and scalability across a range of benchmark tests. The results demonstrate that AMG preconditioners can achieve robust and scalable performance, confirming their potential as practical solvers for large IgA systems in engineering and scientific applications.}

\section{Introduction}
\label{sec:introduction}

In the context of Isogeometric Analysis (IgA), the discretization of a second‐order elliptic boundary‐value problem via tensor‐product B‐splines or NURBS of degree \(p\) on a mesh with \(n\) elements per parametric direction yields a linear algebraic system  
\begin{equation}\label{eq:generic_system}
K \,\mathbf{u} = \mathbf{F},    
\end{equation}
where \(K\in\mathbb{R}^{N\times N}\) is large (with {\(N\propto (n+p)^d\)} in spatial dimension \(d\)), symmetric positive‐definite, and sparse with a bandwidth growing like {$\mathcal{O}( p (n+p)^{d-1})$} %
The improved approximation properties of isogeometric spaces, however, come at the cost of increased coupling between degrees of freedom, which can lead to a higher condition number scaling. As a result, direct solvers may become prohibitively expensive in terms of both memory and arithmetic complexity for three‐dimensional or high‐polynomial‐degree problems. Consequently, Krylov subspace methods~\cite{MR1990645}, combined with specialized preconditioners—such as multigrid techniques adapted to the isogeometric setting~\cite{Manni1,GenericAMG,Mazza1,MR3686804,MR4454927}, overlapping Schwarz methods~\cite{Veiga3,Veiga2,Veiga}, and other domain-decomposition–based approaches~\cite{BPXBuffa,IETIfirst} are essential to achieve scalable and robust solvers.
In particular, preconditioners based on the approximate solution of related matrix equations~\cite{Sangalli2,Sangalli1,Sangalli3,Sangalli4} aim to exploit the tensor‐product structure and the hierarchical nature of spline spaces, thereby ensuring optimal complexity \(\mathcal{O}(N)\) or \(\mathcal{O}(N\log N)\) per iteration. Such techniques enable the efficient simulation of large‐scale elliptic problems in engineering and applied sciences~\cite{VuikStokes}. 

In this work, we aim to assess the effectiveness of Algebraic MultiGrid (AMG) preconditioners for the efficient solution of large-scale, sparse linear systems of the form~\eqref{eq:generic_system}. 
IgA is of particular relevance in this context, as it provides a unified framework for geometry representation and numerical discretization by employing the same basis functions used in Computer-Aided Design to approximate the solution of partial differential equations. 
This tight integration between analysis and geometry eliminates geometric approximation errors, enhances both the accuracy and smoothness of the computed solutions, and enables higher-order continuity across elements—features that are particularly advantageous for problems involving complex geometries or requiring strict geometric fidelity. 

We are particularly interested in evaluating the performance of AMG methods when applied to matrices arising from IgA discretizations. 
Specifically, we consider an AMG formulation in which the aggregation and prolongation operators are constructed directly from the coefficients of the stiffness matrix~$K$, independently of the process by which the matrix was generated. 
This purely algebraic formulation deliberately disregards the geometric information typically available in the IgA setting. 
Although such an approach may appear restrictive, it provides clear benefits in black-box scenarios and in cases where the preconditioner must be assembled directly on the physical domain. 
This stands in contrast to other strategies in which the preconditioner is built from discretization matrices defined in the parametric domain—thereby avoiding the effects of the geometric mapping—and subsequently applied to the full system; see, e.g.~\cite{Manni1} for the multigrid context {based on Toeplitz-like decomposition of the underlying matrices, and~\cite{MR3686804} which is based in turn on a stable splitting of the spline space, while we refer to}~\cite{Sangalli3} for an approach employing matrix equation techniques.

In particular, we focus on AMG techniques that employ aggregation based on the interplay between the principle of compatible relaxation and properties of maximum weight matching~\cite{MR4331965,DAmbra2018BootCMatch}. 
In contrast to classical aggregation schemes~\cite{XuZikatanov,MR1393006}, this approach enables the effective treatment of symmetric and positive definite matrices that do not satisfy the M-matrix property. 
Furthermore, our implementation is designed to exploit both distributed-memory parallelism and modern hardware accelerators, such as GPUs, in order to achieve high scalability and performance on contemporary heterogeneous computing architectures.

To this end, we utilize the software infrastructure provided by the Parallel Sparse Computation Toolkit~\cite{DAMBRA2023100463} (\texttt{PSCToolkit}), which facilitates the scalable construction and manipulation of distributed sparse matrices. The toolkit also supports hybrid computational strategies, enabling seamless integration of CPU-based distributed processing with GPU acceleration through NVIDIA CUDA-enabled devices.

Our investigation encompasses both algorithmic performance and scalability analysis, considering a range of problem sizes and model problems. The overarching goal is to assess the viability of AMG methods as robust, high-performance preconditioners in contemporary high-performance computing (HPC) environments for linear systems obtained from the IgA discretization of elliptic operators.

{The main contributions of this work are as follows. We consider a matching-based AMG preconditioner, based on an aggregation strategy already developed and integrated within the \texttt{PSCToolkit}~\cite{DAMBRA2023100463} infrastructure. A key advantage of this approach is that it does not require the stiffness matrix to satisfy the M-matrix property, thus making it suitable for symmetric positive definite systems arising from high-degree IgA discretizations, where this property is typically violated.
We assess the effectiveness of this method in the context of IgA, leveraging the scalable distributed-memory and GPU-accelerated capabilities of \texttt{PSCToolkit} to demonstrate its practical applicability on modern heterogeneous HPC systems.
Finally, through comprehensive numerical experiments on standard benchmark problems, we evaluate the robustness and scalability of the preconditioner across a range of spline degrees and geometries, including strong scaling results up to 512 MPI tasks. Our findings indicate that AMG methods based on compatible weighted matching provide a viable approach for the solution of large-scale IgA systems.}

Section~\ref{sec:linear_systems} provides a concise overview of the formulation of the linear systems resulting from the isogeometric discretization of elliptic partial differential equations. In Section~\ref{sec:amg_prec}, we present the setup of AMG preconditioners based on parallel aggregation through matching strategies. The integration of these preconditioners within the \texttt{PSCToolkit} is detailed in Section~\ref{sec:psctoolkit}. Section~\ref{sec:numerical_experiments} reports on a series of numerical experiments conducted on standard benchmark problems to assess the performance and scalability of the proposed methods. Finally, Section~\ref{sec:conclusions} summarizes the main findings and outlines prospective avenues for future research.
 
\section{Linear systems from IgA discretizations}
\label{sec:linear_systems}

We consider the Poisson boundary-value problem on a domain $\Omega\subset\mathbb{R}^d$ with boundary $\partial\Omega = \Gamma_D \cup \Gamma_N$ (Dirichlet and Neumann parts).  In strong form the Poisson problem is
\begin{equation}\label{eq:poisson}
    - \nabla^2 u = f \quad \text{in }\Omega,\qquad
u = g_D \;\text{ on }\Gamma_D,\qquad
\frac{\partial u}{\partial \mathbf{n}} = g_N  \;\text{ on }\Gamma_N,
\end{equation}
where $f$ is the source, $g_D,g_N$ are prescribed boundary data and $\mathbf{n}$ is the outward looking normal.  We introduce the Sobolev space $\mathrm{H}^1(\Omega)$ of functions with square-integrable gradient and its subspace $\mathrm{H}^1_0(\Omega)=\{v\in \mathrm{H}^1(\Omega):v|_{\Gamma_D}=0\}$.  The weak form is: seek $u\in \mathrm{H}^1(\Omega)$ satisfying $u|_{\Gamma_D}=g_D$ such that
\[
\int_\Omega \nabla u\cdot\nabla v \,\mathrm{d}\Omega
= \int_\Omega f\,v\,\mathrm{d}\Omega + \int_{\Gamma_N} g_N\,v\,\mathrm{d}\Gamma
\quad\forall v\in \mathrm{H}^1_0(\Omega).
\]
Equivalently, one can write $a(u,v)=L(v)$ with 
$a(u,v) = \int_\Omega \nabla u\cdot\nabla v\,\mathrm{d}\Omega$ and 
$L(v)=\int_\Omega f\,v\,\mathrm{d}\Omega + \int_{\Gamma_N} g_N\,v\,\mathrm{d}\Gamma$.  This is the standard weak formulation of the Poisson problem.  The solution of this variational problem (under suitable conditions on the data) is the weak solution of the original PDE; see, e.g.~\cite[\S 6]{Evans}.

\subsection{Isogeometric preliminaries: function spaces, bases, and mappings}
{All constructions in this section are presented at the level of a single parametric patch; the extension to multi-patch domains will be discussed in Section~\ref{sec:multipatch}.}

Isogeometric analysis (IgA) uses spline or NURBS basis functions defined on a reference domain and a mapping to the physical domain.  Let $\widehat\Omega=[0,1]^d$ be the parametric domain.  A (tensor-product) B-spline basis $\{N_i(\boldsymbol{\xi})\}$ or NURBS basis of degree $p$ is defined by knot vectors in each coordinate~\cite[\S 12]{Schumaker}. In particular, a NURBS basis function~\cite[\S 4]{zbMATH00792231} is 
\[
R_i(\boldsymbol{\xi}) = \frac{w_i N_i(\boldsymbol{\xi})}{\sum_j w_j N_j(\boldsymbol{\xi})},
\]
where the $w_i>0$ are weights.  The geometry mapping is given by 
\[
\boldsymbol{x} = F(\boldsymbol{\xi}) = \sum_{i} \mathbf{P}_i\,R_i(\boldsymbol{\xi}),
\]
with control points $\{\mathbf{P}_i\}$.  

{In IgA one often adopts the isoparametric paradigm, in which the same basis is used to describe both the geometry and the discrete fields. On a given patch, the associated physical basis functions are obtained by composition with the geometry mapping, i.e.,
\[
\phi_i(\boldsymbol{x}) = R_i(\boldsymbol{\xi}), \qquad \boldsymbol{x}=F(\boldsymbol{\xi}).
\]
These mapped functions inherit the partition-of-unity and local-support properties of the underlying splines.}

{In a completely analogous manner, the entire construction can be carried out using a B-spline basis $\{N_i(\boldsymbol{\xi})\}$, without introducing weights. In this case, both the geometry mapping and the discrete space are defined directly in terms of B-splines, rather than rational functions. Which is not the special case $w_i=1$ in the NURBS definition, but corresponds to a formulation in which no rational representation is employed.}

\subsection{Single--patch IgA discretization}
{We now consider a single patch domain $\Omega = F(\widehat\Omega)$. Using the basis functions introduced above, we define the discrete space $S_h \subset \mathrm{H}^1(\Omega)$ as the span of the physical basis functions $\{\phi_i\}$.}

Applying the Galerkin method ($u_h=\sum_i u_i\phi_i$) yields a linear system $K\mathbf{u}=\mathbf{F}$ with entries
\[
k_{ij} = \int_\Omega \nabla \phi_j\cdot\nabla \phi_i \,\mathrm{d}\mathbf{x}, 
\qquad
f_i = \int_\Omega f\,\phi_i\,\mathrm{d}\mathbf{x} + \int_{\Gamma_N} g_N\,\phi_i\,\mathrm{d}S.
\]
By change of variables this is evaluated on $\widehat\Omega$.  In particular
\[
k_{ij} 
= \int_{\widehat\Omega} (\nabla_{\boldsymbol{\xi}}R_j(\boldsymbol{\xi}))^\top 
\bigl(J_F^{-T}J_F^{-1}\bigr)\,\nabla_{\boldsymbol{\xi}}R_i(\boldsymbol{\xi}) \,\det(J_F)\,\mathrm{d}\boldsymbol{\xi},
\]
where $J_F=\partial F/\partial\boldsymbol{\xi}$ is the Jacobian of the geometry mapping.  The load vector entries are similarly computed by pulling back the integrals of $f\phi_i$ and $g_N\phi_i$ to the parameter domain. The matrix assembly follows the standard Galerkin procedure: one enumerates the basis functions, computes element-level integrals, and assembles into the global matrix.  Concretely, for each knot-span element $e\subset\widehat\Omega$, we compute 
  $$k^e_{mn} = \int_e \nabla\phi_m\cdot\nabla\phi_n \,\mathrm{d}\Omega$$ and 
  $$f^e_m = \int_e f\,\phi_m \,\mathrm{d}\Omega + \int_{\partial e\cap\Gamma_N} g_N\,\phi_m \,\mathrm{d}S$$ 
  using quadrature (pulling back via $F$). Then, we add the local contributions $k^e_{mn},f^e_m$ into the global $k_{ij},f_i$ using the mapping from local to global basis indices. For Dirichlet boundary conditions we identify basis functions whose support lies on $\Gamma_D$ (e.g.\ control points on the boundary) and enforce $u_i=g_D$ there.  In practice one fixes those DoFs by modifying $K$ and $F$ (e.g.\ row/column elimination or penalty). On the other hand, Neumann BCs are incorporated via the boundary integrals in $f_i$ on $\Gamma_N$.
The result is a symmetric positive-definite system. Note that NURBS and B-Spline basis are globally smooth but locally supported: each basis spans multiple elements, so $K$ has a sparsity pattern similar to a high-order finite-element matrix; see the left panel in Fig.~\ref{fig:pattern}.

{A straightforward application of the standard FEM assembly loop is computationally inefficient and entails a significant memory footprint, due to the use of large local matrices and high-order quadrature rules. More efficient alternatives are available, including generalized Gaussian quadrature~\cite{Quad1}, sum-factorization techniques exploiting the tensor-product structure of splines~\cite{Quad2,MR4087175,MR3951499}, and their combinations, which can substantially reduce the computational cost and memory requirements~\cite{Quad3}.}

\subsection{Multi--patch IgA discretization}\label{sec:multipatch}
For complex geometries, the domain $\Omega$ is covered by multiple patches $\{\Omega^{(p)}\}$.  Each patch $p$ has its own mapping $F^{(p)}:\widehat\Omega\to\Omega^{(p)}$ and local basis $\{R^{(p)}_i\}$.  On patch $p$ we define the physical basis functions $\phi^{(p)}_i(\boldsymbol{x}) = R^{(p)}_i(\boldsymbol{\xi})$ with $\boldsymbol{x}=F^{(p)}(\boldsymbol{\xi})$.  The global function space is constructed by “gluing” together the individual patch spaces.
In a conforming multi-patch design, the parameterizations along shared boundaries are aligned and the basis functions on each interface are correspondingly identified.

\begin{figure}[htb]
\sidecaption
\includegraphics[width=3.6cm]{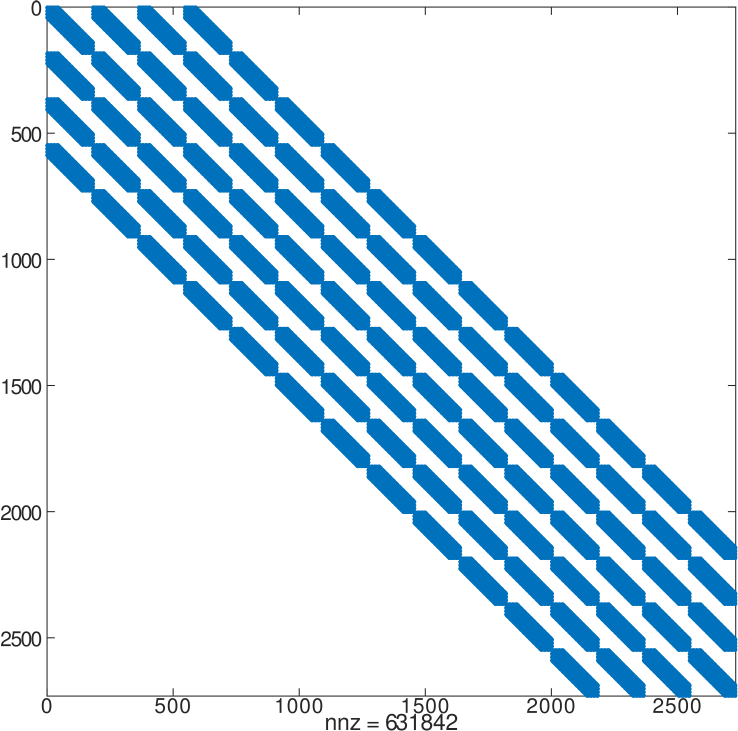}
\includegraphics[width=3.6cm]{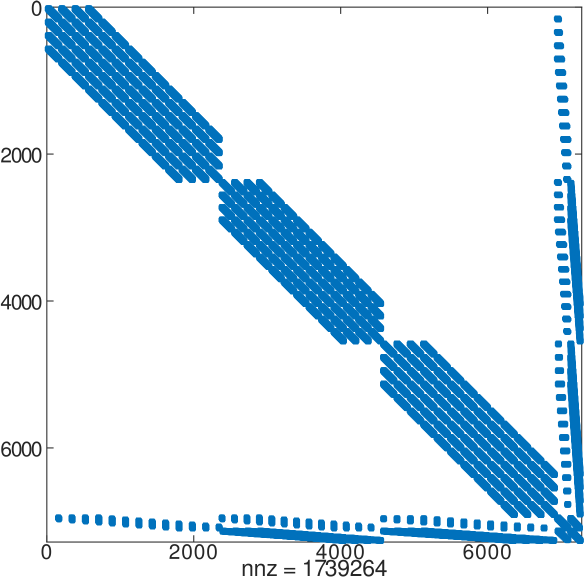}
\caption{Spy patterns of an IgA stiffnes matrix for the Poisson equation on a cube with B-Spine basis of degree 3 and regularity 2, and of the Poisson equation on a thick $L$-shaped domain---with the same degree and regularity---but built with three patches.}
\label{fig:pattern}       %
\end{figure}

Assembly proceeds patch-by-patch, analogous to multi-element FEM, we repeat the procedure for the single patch on each local patch $p$, computing the local stiffness matrix, i.e., on each patch  
$$k^{(p)}_{ij} = \int_{\Omega^{(p)}} \nabla \phi^{(p)}_j \cdot \nabla \phi^{(p)}_i \,\mathrm{d}\mathbf{x}$$ 
and load vector 
$$g^{(p)}_i = \int_{\Omega^{(p)}} f\,\phi^{(p)}_i\,\mathrm{d}\mathbf{x} + \int_{\partial\Omega^{(p)}\cap\Gamma_N} g_N\,\phi^{(p)}_i\,\mathrm{d}S$$
by pulling back to the parametric domain of patch $p$. Then a global assembly starts, and we add the contributions $k^{(p)}_{ij}$ to the global matrix, mapping local patch indices to global DoFs.  When patches share an interface, the basis functions on that interface from each patch are identified as the same global DoF. At last, Dirichlet boundary values are imposed on the basis functions supported on the external boundary $\Gamma_D$, while Neumann data on $\Gamma_N$ are incorporated through boundary integrals over each patch. The resulting system $K\mathbf{u}=\mathbf{F}$ couples all patches together. Its sparsity pattern reflects patch connectivity: blocks corresponding to each patch, with coupling at shared interfaces; see the right panel in Fig.~\ref{fig:pattern} in which the blocking structure plus interfaces is clearly highlighted.
{We note that the construction of the AMG preconditioner, discussed in Section~\ref{sec:amg_prec}, is purely algebraic; therefore, the ordering of the DoFs, particularly in the multipatch case, is not exploited in the solution procedure.}

\subsection{Isoparametric vs.\ non--isoparametric IgA}\label{sec:iso-or-noniso}
{By default, IgA adopts the isoparametric paradigm, in which the geometry mapping and the discrete solution are represented within the same spline space defined on the reference domain. In particular, the geometry is assumed to be a function in the span of the chosen spline basis, and the same space is used to construct the discrete trial and test functions.}
In this case the assembly is as above.  One uses the Jacobian of the geometric mapping and the same basis gradients for the stiffness integrals.  For example, on a patch the entries are 
\[
k_{ij} = \int_{\Omega^{(p)}} \nabla R_j\cdot\nabla R_i \,\mathrm{d}\mathbf{x},
\]
with $R_i$ the NURBS/B-Spline basis from the geometry, just as in the single-patch formula.  There is no distinction between geometry and solution bases in the assembly~\cite{HUGHES2021467}.

A {\it non-isoparametric} discretization would use different basis functions for the solution than those defining the geometry.  For instance, the geometry might use NURBS $\{R_i\}$ while the solution uses a distinct B-Spline basis $\{\widetilde N_i\}$.  In that case $\phi_i(\boldsymbol{x})=\widetilde N_i(\boldsymbol{\xi})$ and the geometry mapping $F$ (and its Jacobian $J_F$) are determined by $\{R_i\}$~\cite{HUGHES2021467}.  Non-isoparametric isogeometric spaces have been extensively and successfully employed in the context of structure-preserving discretizations, particularly for the approximation of vector- and tensor-valued fields defined on manifolds (see, e.g.~\cite{Arnold20061,MR2808112,Hiemstra20141444,Evans2020422}). Within these frameworks, scalar and vector fields are interpreted as discrete analogues of differential forms, and their push-forwards from the parametric to the physical domains are constructed in a manner consistent with their geometric character, hence it makes sense to include the non-isoparametric setting into the analysis. The stiffness entries become
\[
k_{ij} = \int_{\widehat\Omega} (\nabla_{\boldsymbol{\xi}}\widetilde N_j)^\top (J_F^{-T}J_F^{-1})\,(\nabla_{\boldsymbol{\xi}}\widetilde N_i)\,\det(J_F)\,\mathrm{d}\boldsymbol{\xi},
\]
which differs in that the gradient is of $\widetilde N_i$ instead of $R_i$.  %
{In practice, isogeometric formulations often adopt an isoparametric approach, in which the same spline space is used for both the geometry and the discrete solution. However, this choice is not universal. In the literature, one can find numerous examples where different spline spaces are employed~\cite{Arnold20061,MR2808112,Hiemstra20141444,Evans2020422}, such as NURBS for the geometry and B-splines for the solution, or spaces with differing polynomial degrees and smoothness. The isoparametric setting is therefore a common, but not mandatory, choice.}
{Boundary conditions are enforced in the same way. For homogeneous Dirichlet conditions, the corresponding degrees of freedom associated with basis functions supported on $\Gamma_D$ are fixed to zero. In the case of non-homogeneous Dirichlet data, standard approaches such as a lifting of the boundary data or a modification of the right-hand side are employed to incorporate the prescribed values; see~\cite[Section~3.7.2]{GeoPDEs2}.}

In summary, the two approaches differ only in the choice of basis; in the isoparametric IgA, geometry and solution use the same basis; this permits using exact CAD geometry and high smoothness while performing the assembly by means of a single Jacobian. In the non-isoparametric IgA formulation, geometry and solution bases differ.  One still computes $k_{ij}=\int_\Omega\nabla\phi_j\cdot\nabla\phi_i$, but with $\phi_i$ from the alternative basis while $J_F$ from the geometry mapping.

\section{AMG preconditioners}
\label{sec:amg_prec}

In all cases discussed in Section~\ref{sec:linear_systems}, the IgA discretization leads to a Galerkin system of the form~\eqref{eq:generic_system} with
\begin{equation}\label{eq:matrix_and_vector_values}
  k_{ij} = \int_\Omega \nabla \phi_j \cdot \nabla \phi_i \,\mathrm{d}\mathbf{x},
  \qquad
  f_i = \int_\Omega f\,\phi_i\,\mathrm{d}\mathbf{x} + \int_{\Gamma_N} g_N\,\phi_i\,\mathrm{d}S,  
\end{equation}
where $\{\phi_i\}$ denotes the global IgA basis on $\Omega$.  
The basis functions $\phi_i$ and the mapping $F$ are chosen according to the single- or multi-patch geometry and the isoparametric or non-isoparametric formulation.  
{Dirichlet boundary conditions are imposed by constraining the degrees of freedom associated with basis functions supported on $\Gamma_D$ as described in~\cite[Section~3.7.2]{GeoPDEs2}.}
Aside from these modeling choices, the assembly of $K$ and $\mathbf{F}$ follows the same variational procedure in each case and shares the property of being symmetric, positive definite, and sparse.  

In this configuration, a natural choice is to adopt a Krylov method, specifically the Conjugate Gradient (CG) method~\cite[\S 6.7]{MR1990645}, coupled with a suitable preconditioner to accelerate convergence, typically of the multigrid type~\cite{XuZikatanov}.  
For our $K \in \mathbb{R}^{N\times N}$, this consists in considering a nested sequence of aggregation levels indexed by $\ell=0,1,\dots,L$, with level $0$ the finest and level $L$ the coarsest.  
On each level $\ell$ we have a matrix $K^{(\ell)}\in\mathbb{R}^{N_\ell\times N_\ell}$, with $N_0=N$ and $N_L\ll N$, and introduce prolongation operators 
\[
P^{(\ell+1)}_{\ell}\colon\mathbb{R}^{N_{\ell+1}}\to\mathbb{R}^{N_{\ell}}
\]
and corresponding restriction operators {$R^{(\ell)}_{\ell+1}=(P^{(\ell+1)}_{\ell})^\top$.}
These are chosen so that the Galerkin coarse‐grid operator  
\begin{equation}\label{eq:galerkin_matrix_projection}
  K^{(\ell+1)} \;=\; R^{(\ell)}_{\ell+1}\,K^{(\ell)}\,P^{(\ell+1)}_{\ell}
\end{equation}
remains symmetric positive definite.
The pair $(P,R)$ transfers residuals and corrections between levels, ensuring compatibility of energy norms across the hierarchy.  

On each level $\ell$, a smoother approximating the action of $(K^{(\ell)})^{-1}$ damps high‐frequency error components cheaply, i.e., with linear cost in the number of unknowns on that level and limited interprocessor communication.  
Denote by $M^{(\ell)}_{\text{pre}}$ and $M^{(\ell)}_{\text{post}}$ the pre‐ and post‐smoothing operators, respectively.  
A single smoothing sweep updates
\begin{equation}
\label{smoother}
   \mathbf u^{(\ell)} \;\leftarrow\; \mathbf u^{(\ell)} + M^{(\ell)}\bigl(\mathbf f^{(\ell)} - K^{(\ell)}\mathbf u^{(\ell)}\bigr),
\end{equation}
where $M^{(\ell)}$ is either $M^{(\ell)}_{\text{pre}}$ or $M^{(\ell)}_{\text{post}}$, depending on the stage of the cycle.

{
A two-level correction then consists of:
\begin{enumerate}
  \item \textbf{Pre-smoothing:} perform $\nu_1$ iterations of scheme~\eqref{smoother} applied to $(K^{(\ell)},\mathbf f^{(\ell)})$, using~$M^{(\ell)}_{\text{pre}}$;
  \item \textbf{Restrict residual:} $\mathbf r^{(\ell+1)} = R^{(\ell)}_{\ell+1}\bigl(\mathbf f^{(\ell)} - K^{(\ell)}\mathbf u^{(\ell)}\bigr)$;
  \item \textbf{Coarse solve:} approximately solve $K^{(\ell+1)}\mathbf e^{(\ell+1)} = \mathbf r^{(\ell+1)}$. This system is solved either by a direct method or by an iterative solver;
  \item \textbf{Prolongate and correct:} $\mathbf u^{(\ell)} \leftarrow \mathbf u^{(\ell)} + P^{(\ell+1)}_{\ell}\,\mathbf e^{(\ell+1)}$;
  \item \textbf{Post-smoothing:} perform $\nu_2$ iterations of scheme~~\eqref{smoother} applied to $(K^{(\ell)},\mathbf f^{(\ell)})$, using~$M^{(\ell)}_{\text{post}}$.
\end{enumerate}}

Embedding this two‐level idea recursively yields the V‐cycle operator $\mathcal{P}^{-1}$.  
The resulting operator $\mathcal{P}^{-1}\approx K^{-1}$ is symmetric positive definite when the restriction and prolongation operators satisfy $R=P^\top$ and the post‐smoother is chosen as the adjoint (or transpose) of the pre‐smoother.  
In this case, the V‐cycle defines a symmetric iteration suitable for use as a preconditioner in the CG method applied to $K\mathbf u=\mathbf F$; see~\cite{XuZikatanov}.

The theoretical objective is for the preconditioned system  
\begin{equation}
   \mathcal{P}^{-1}K\,\mathbf u = \mathcal{P}^{-1}\mathbf F
\end{equation}
to have a condition number bounded independently of $N$, ensuring \emph{optimal convergence} of the iterative solver (i.e., a convergence rate independent of problem size). {In addition, for isogeometric discretizations, robustness is also expected with respect to geometric parametrization $F$, spline degree $p$, and inter-element smoothness, so that the condition number remains bounded independently not only of $N$, but also of these discretization-specific parameters; see the discussion on this subject in the case of geometric multigrid methods in~\cite{Manni1,MR4454927,MR3686804}.} Optimal convergence is distinct from \emph{optimal computational complexity}, which refers to the total cost per iteration scaling linearly with $N$.  
Both aspects are desirable but not equivalent, and practical AMG implementations often achieve near‐optimal performance only under additional assumptions on the smoother and the coarsening strategy.

From an implementation standpoint, all operations composing the method should be suitable for parallel execution, consisting mostly of local computations to minimize communication in distributed‐memory environments, and capable of exploiting Single‐Instruction–Multiple‐Data (SIMD) parallelism offered by multicore processors and GPU accelerators.

In our case, we focus on constructing an AMG method in which the aggregation and prolongation operators are derived directly from the coefficients of the matrix \( K \), regardless of how the matrix was generated.  
This approach deliberately ignores geometric information typically available in the IgA context.  
While this may seem limiting, it is advantageous in black‐box scenarios or when constructing the preconditioner directly on the physical domain.  
In contrast, the more common approach builds the preconditioner on the discretization matrices defined in the parametric domain—thus avoiding the influence of the geometric mapping—and subsequently applies it to the full system; see, e.g.~\cite{Manni1,Sangalli3}.

\subsection{Aggregation with compatible weighted matching}
\label{sec:match_aggreg}

In the \emph{compatible relaxation} (CR) framework~\cite{XuZikatanov},
the quality of a coarse space is assessed by performing a relaxation that acts only on
the \emph{fine} variables while keeping the \emph{coarse} ones invariant.
If this restricted relaxation converges rapidly, the chosen coarse space is said to be
\emph{compatible}, meaning that it captures the slowly varying components of the error.
This principle motivates the aggregation strategy known as
\emph{compatible weighted matching}, first
implemented in the BootCMatch software framework~\cite{DAmbra2018BootCMatch}, for which a parallel version is available in PSCToolkit.

\paragraph{Edge weights and projection onto the fine subspace.}
Let \(K\in\mathbb{R}^{N\times N}\) be a symmetric positive definite matrix and
\(G=(\mathcal V,\mathcal E)\) its adjacency graph,
where vertices correspond to unknowns and edges to nonzero off-diagonal entries of \(K\).

{
Given a \emph{test vector} $\mathbf w \in \mathbb{R}^N$ approximating an algebraically smooth error, the goal is to group together variables whose components are similar in $\mathbf w$.
The test vector $\mathbf w$ is obtained by applying a few iterations of a relaxation scheme (e.g., damped Jacobi or Gauss--Seidel) to the homogeneous system
\[
K \mathbf w = 0,
\]
starting from a random initial vector. This process damps the high-frequency components of the error, so that $\mathbf w$ represents an approximation of the algebraically smooth error components.
We note that the test vector is computed on the finest level during the setup phase and is then recursively restricted to the coarser levels throughout the AMG hierarchy construction.
}

For each edge \((i,j)\in\mathcal E\), the compatibility coefficient is defined as
\begin{equation}\label{eq:cij}
  c_{ij}
  = 1 - \frac{2\,k_{ij} w_i w_j}{k_{ii} w_i^2 + k_{jj} w_j^2},
  \qquad (i,j)\in\mathcal E .
\end{equation}
This coefficient originates from the analysis of the \emph{fine-space projection}
associated with the local two-variable subsystem
\[
  K_e =
  \begin{pmatrix}
    k_{ii} & k_{ij}\\[2pt]
    k_{ij} & k_{jj}
  \end{pmatrix},
  \qquad
  \mathbf w_e = (w_i, w_j)^\top .
\]

{In our compatible relaxation framework, the coarse space is locally represented by the vector $\mathbf w_e = (w_i, w_j)^\top$.
The associated fine space is defined as the one-dimensional subspace generated by a vector $\mathbf w_e^{\perp}$ that is orthogonal to $\mathbf w_e$ with respect to the energy inner product induced by $K_e$, i.e.,
\[
(\mathbf w_e^{\perp})^\top K_e \mathbf w_e = 0.
\]
Therefore, the compatibility coefficient $c_{ij}$ in~\eqref{eq:cij} measures how well the local behavior of the pair $(i,j)$ is captured by the coarse space spanned by $\mathbf w_e$. Smaller values indicate that the pair $(i,j)$ is well represented by the coarse space (i.e., both variables belong to the same algebraically smooth mode), while larger values indicate that the pair exhibits fine-scale behavior.
}

\paragraph{Maximum product matching and CR convergence.}

{
A \emph{matching} $\mathcal M' \subset \mathcal E$ is a set of edges such that no two edges share a common vertex (Fig.~\ref{fig:matching_example}), i.e., each vertex belongs to at most one edge in $\mathcal M'$~\cite[p.35]{MR3644391}.}

The aggregation is obtained by computing a maximum product matching { (Fig.~\ref{fig:matching_example_max})}:
\begin{equation}\label{eq:maxprod}
  \mathcal M
  = \argmax_{\mathcal M' \subset \mathcal E \text{ matching}}
  \;\prod_{(i,j)\in\mathcal M'} c_{ij}.
\end{equation}
\begin{figure}[htbp]
\centering
\begin{subfigure}{0.3\textwidth}
\centering
\begin{tikzpicture}
\Vertex[x=0,y=1,size=0.2,color=none]{A}
\Vertex[x=1,y=2,size=0.2,color=none]{B}
\Vertex[x=2,y=1,size=0.2,color=none]{C}
\Vertex[x=1,y=0,size=0.2,color=none]{D}
\Edge[label=5](A)(B)
\Edge[label=3](B)(C)
\Edge[label=4](C)(D)
\Edge[label=2](D)(A)
\Edge[label=6](A)(C)
\Edge[label=1](B)(D)
\end{tikzpicture}
\caption{Weighted $\mathcal{G}=(\mathcal{V},\mathcal{E})$}
\end{subfigure}
\begin{subfigure}{0.3\textwidth}
\centering
\begin{tikzpicture}
\Vertex[x=0,y=1,size=0.2,color=none]{A}
\Vertex[x=1,y=2,size=0.2,color=none]{B}
\Vertex[x=2,y=1,size=0.2,color=none]{C}
\Vertex[x=1,y=0,size=0.2,color=none]{D}
\Edge[color=gray!30](A)(B)
\Edge[color=gray!30](B)(C)
\Edge[color=gray!30](C)(D)
\Edge[color=gray!30](D)(A)
\Edge[color=gray!30](A)(C)
\Edge[color=gray!30](B)(D)
\Edge[lw=2pt,color=blue](A)(C)
\Edge[lw=2pt,color=blue](B)(D)
\end{tikzpicture}
\caption{$\mathcal{M}' \subset \mathcal{E}$\label{fig:matching_example}}
\end{subfigure}
\begin{subfigure}{0.3\textwidth}
\centering
\begin{tikzpicture}
\Vertex[x=0,y=1,size=0.2,color=none]{A}
\Vertex[x=1,y=2,size=0.2,color=none]{B}
\Vertex[x=2,y=1,size=0.2,color=none]{C}
\Vertex[x=1,y=0,size=0.2,color=none]{D}
\Edge[color=gray!30](A)(B)
\Edge[color=gray!30](B)(C)
\Edge[color=gray!30](C)(D)
\Edge[color=gray!30](D)(A)
\Edge[color=gray!30](A)(C)
\Edge[color=gray!30](B)(D)
\Edge[lw=2pt,color=red](A)(B)
\Edge[lw=2pt,color=red](C)(D)
\end{tikzpicture}
\caption{$\mathcal{M}$}
\end{subfigure}
\caption{{(a) Graph, (b) a matching $\mathcal{M}' \subset \mathcal{E}$ (product of the edge weights is $1$), (c) a maximum product matching $\mathcal{M}$ (product of the edge weights is $20$.)}\label{fig:matching_example_max}}
\end{figure}

Maximizing the product of the compatibility coefficients identifies a set of disjoint
edges that define a coarse space whose complementary fine space, constructed as the
orthogonal complement (in the energy inner product induced by \(K\)) of the coarse
space, is well conditioned, then a
relaxation applied only to the fine variables while keeping the coarse ones invariant
converges rapidly.
This is the fundamental rationale of the {\em compatible weighted matching} approach: the
maximum product matching ensures an energetically orthogonal decomposition of the space
into coarse and fine components, such that the fine-space relaxation is effective.

{
In the practice, the maximum product matching is computed by applying a linear-time approximate graph matching algorithm~\cite{6009071} (e.g., a greedy maximum weight matching algorithm applied to the logarithms of the weights).
}

\paragraph{Interpolation and coarse operator.}
Each matched edge \((i,j)\in\mathcal M\) defines a two-point aggregate
\(\mathcal G_p=\{i,j\}\), while any possible unmatched vertex (related to approximate maximum matching) forms a singleton aggregate.
After reordering the unknowns so that singleton indices appear last,
let \(n_p = |\mathcal M|\), \(n_s\) be the number of singletons,
and \(n_c = n_p + n_s\).
The prolongation operator is constructed as
\begin{equation}\label{eq:prolongation}
P
=
\begin{pmatrix}
\operatorname{blockdiag}(\mathbf{w}_{e_1},\dots,\mathbf{w}_{e_{n_p}}) & 0\\[6pt]
0 & \operatorname{diag}(w_s/|w_s|)
\end{pmatrix}
\in\mathbb{R}^{N\times n_c},
\quad
\mathbf{w}_e = \frac{(w_i,w_j)^\top}{\|(w_i,w_j)\|_2},
\end{equation}
and the corresponding coarse matrix is obtained by the Galerkin projection
\begin{equation*}
  K_c = P^\top K P.
\end{equation*}
Recursively applying the same procedure to \(K_c\) yields a multilevel AMG hierarchy.
Note that the reordering in \eqref{eq:prolongation} is introduced solely for explanatory clarity; no such reordering is performed in the actual implementation.

\paragraph{Variants and efficiency.}
The same framework allows $m$–step aggregation by replacing \(P\) with
\(P_1P_2\cdots P_m\), thus forming aggregates of size up to \(2^m\).
Convergence can be further improved by applying a smoothing step to the interpolation,
\(\overline P=(I-\omega D^{-1}K)P\),
or by applying the hierarchy within a K–cycle accelerated by a CG iteration,
and then using it as a preconditioner inside the Flexible CG method~\cite{FCG}.
The efficiency of the resulting hierarchy is characterized by the
\emph{operator complexity}
\begin{equation}
\operatorname{OPC}
= \frac{\sum_{\ell=0}^{L} \operatorname{nnz}(K^{(\ell)})}
       {\operatorname{nnz}(K^{(0)})},
\end{equation}
where $\operatorname{nnz}(K^{(\ell)})$ denotes the number of nonzero entries
in the system matrix~\eqref{eq:galerkin_matrix_projection} at level~$\ell$,
and $L$ is the total number of levels.
A lower $\operatorname{OPC}$ indicates a cheaper hierarchy,
while higher values generally lead to more accurate coarse-grid corrections.

\begin{figure}
\sidecaption
\includegraphics[width=7.8cm]{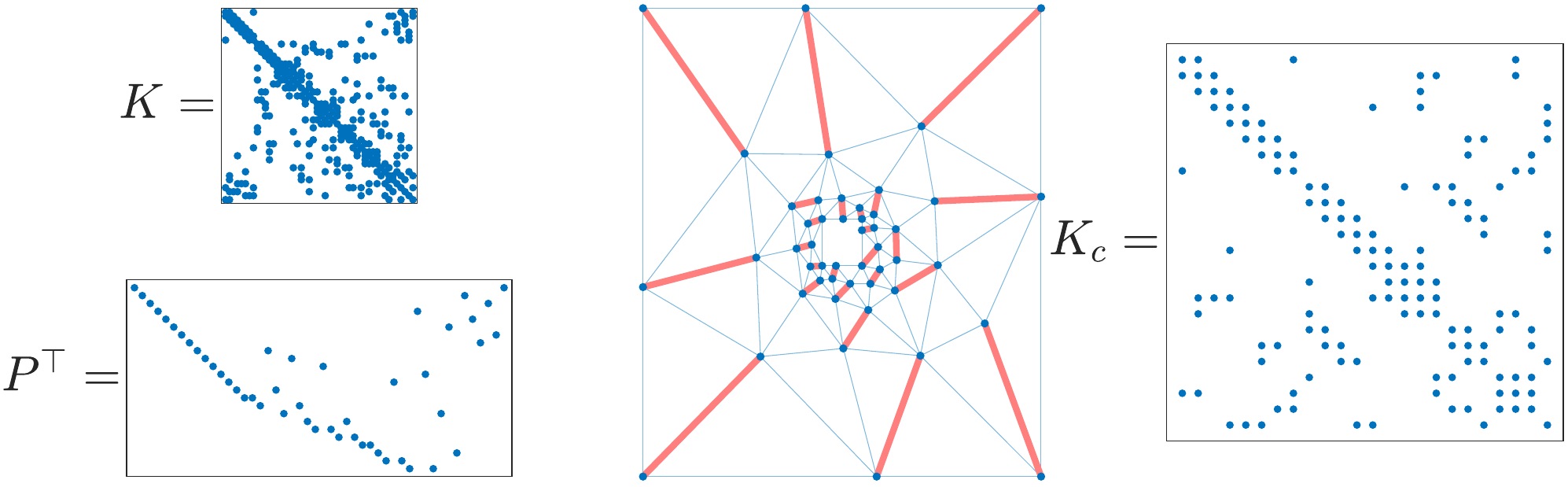}
\caption{Example of compatible weighted graph matching aggregation:
from the weighted adjacency graph (edge weights $c_{ij}$),
a maximum product matching (bold edges) is computed,
the interpolation operator \(P\) is assembled,
and the coarse matrix is obtained as \(K_c = P^\top K P\).
The resulting $K-$orthogonal operator (i.e., the fine space operator) is more diagonally dominant,
ensuring rapid convergence of the compatible relaxation.}
\label{fig:matching_aggregation}
\end{figure}

\subsection{Smoothers via polynomial acceleration of stationary iterative solvers}
The remaining key ingredient in designing an efficient AMG solver is the choice of a suitable smoother.  In our parallel setting, it is crucial that the smoother’s action remain largely local or be implementable solely via sparse‐matrix–vector products (\texttt{spmv}) and vector updates (\texttt{axpby})~\cite{MR2861652}.  For a purely distributed implementation---without relying on multithreading or GPU acceleration---a natural option is a block‐Jacobi smoother in which each MPI rank inverts its diagonal block of DoFs via a fixed‐point iteration (e.g., Gauss-Seidel).  While this block‐Jacobi plus Gauss–Seidel combination is robust and easy to integrate, solving sparse triangular systems can become a bottleneck when targeting SIMD architectures.  

To alleviate this, one can employ fully parallel smoothers such as the $\ell_1$‐Jacobi method, which replaces exact block inversions by inexpensive diagonal scaling and offloads all work to embarrassingly parallel \texttt{spmv}/\texttt{axpby} kernels.  Alternatively, approximate‐inverse techniques---where each block’s inverse is precomputed in sparse, factorized form---offer a compromise between convergence speed and vectorization efficiency.  Polynomial smoothers~\cite{Lottes} (e.g., Chebyshev relaxations) also fit naturally into this framework: they require only matrix–vector products and recurrence coefficients chosen considering the optimization of the AMG convergence bounds.   

The choice among these smoothers should be guided by the target hardware’s arithmetic intensity, memory‐bandwidth characteristics, and the spectral properties of the matrix---tuning relaxation weights or polynomial degrees to balance convergence rates against per‐iteration cost.  Such flexibility ensures that the AMG preconditioner remains both scalable and efficient across a wide range of distributed‐memory platforms.

In here we focus on the usage of the \(\ell_1\)\nobreakdash-Jacobi smoother~\cite{MR2861652} which is, as we mentioned, a fully parallel relaxation scheme designed to avoid any triangular solves by replacing each diagonal entry \(k_{ii}\) of the matrix \(K\) with an \(\ell_1\)-norm‐based scaling factor.  Specifically, one defines
\[
d_{i} \;=\; k_{ii} \;+\; \sum_{j\neq i} \lvert k_{ij}\bigr \rvert, \quad D = \operatorname{diag}(d_{i}),
\]
and then performs the fixed-point iteration
\[
\mathbf{x}^{(k+1)} \;=\; \mathbf{x}^{(k)} \;+\; D^{-1}(\mathbf{f} - K\,\mathbf{x}^{(k)}),
\]
Because \(D\) contains only diagonal entries, each update reduces to one sparse‐matrix–vector multiply (\texttt{spmv}) to compute \(K\,\mathbf{x}^{(k)}\) and two vector \texttt{axpby}‐type operations.  The inclusion of the off-diagonal absolute sums in \(D_{ii}\) provides stronger damping of high‐frequency error components than classical Jacobi, often yielding faster convergence in practice~\cite{MR2861652,MR4331965}, yet retains perfect data locality and vectorization potential on modern SIMD architectures. To improve the robustness of the Jacobi iteration, we employ the Chebyshev acceleration from~\cite{DAmbra2025}, which requires no additional spectral estimates of $K$ and adds no extra sparse‐matrix–vector products---incurring only one additional \texttt{axpby} operation per smoother sweep, and are an extension of the smoothers based on 4\textsuperscript{th}-kind Chebyhsev acceleration introduced in~~\cite{Lottes}.

\section{PSCToolkit implementation}
\label{sec:psctoolkit}

The Parallel Sparse Computation Toolkit\footnote{Documentation and installation instructions are available at \href{https://psctoolkit.github.io}{psctoolkit.github.io}.}~\cite{DAMBRA2023100463} (\texttt{PSCToolkit}) is a software framework for large-scale sparse computations in distributed-memory environments.  
It is written mainly in modern Fortran~2008 to ensure code stability and performance portability across HPC architectures.
The toolkit consists of two modular libraries:
\begin{description}[PSBLAS]
  \item[\texttt{PSBLAS},]implementing distributed sparse-BLAS operations;
  \item[\texttt{AMG4PSBLAS},]providing algebraic multigrid (AMG) preconditioners and solvers.
\end{description}
MPI handles inter-process communication, while OpenMP and CUDA provide shared-memory and GPU parallelism.  
The build system supports both \texttt{GNU Automake} and \texttt{CMake}, automatically detecting backends (MPI, OpenMP, CUDA, and optionally OpenACC) and setting compile-time flags for either a full or per-module build.

\paragraph{Distributed data structures.}
Parallel execution is managed through a \lstinline[language=Fortran,style=fortran]{type(psb_ctxt_type)} context wrapping an MPI communicator.  
The degrees of freedom (DoFs) are distributed using a descriptor:
\begin{lstlisting}[language=Fortran,style=fortran]
type(psb_desc_type) :: desc
call psb_cdall(ctxt, desc, info, vl=vl)
\end{lstlisting}
where \lstinline[language=Fortran,style=fortran]{vl} contains the local list of global indices, e.g., from a METIS partition~\cite{METIS}.  
A distributed sparse matrix is then allocated and filled as:
\begin{lstlisting}[language=Fortran,style=fortran]
type(psb_dspmat_type) :: K
call psb_spall(K, desc, info, nnz=nnz)
call psb_spins(num_coeffs, irow, icol, val, K, desc, info)
\end{lstlisting}
Overlapping entries contributed by multiple processes are automatically summed during assembly:
\[
K_{ij} = \sum_{\ell=1}^{r} k_{ij}^{p_\ell}.
\]
A similar workflow applies to right-hand-side vectors:
\begin{lstlisting}[language=Fortran,style=fortran]
type(psb_d_vect_type) :: f
call psb_geall(f, desc, info)
call psb_geins(num_coeffs, irow, val, f, desc, info)
\end{lstlisting}

\paragraph{Finite element assembly.}
In a finite element context, each process loops over its local elements, computes local stiffness matrices \(\mathbf{k}^e\) and load vectors \(\mathbf{f}^e\), and inserts their contributions into the global structures using \texttt{psb\_spins} and \texttt{psb\_geins}.  
Global consistency is ensured by three collective calls:
\begin{lstlisting}[language=Fortran,style=fortran]
call psb_cdasb(desc, info)
call psb_spasb(K, desc, info)
call psb_geasb(f, desc, info)
\end{lstlisting}
After these steps, the fully assembled system \( K\mathbf{u} = \mathbf{F} \) is ready for solution.

\paragraph{GPU backends.}
If no \lstinline[language=Fortran,style=fortran]{mold} in the \lstinline[language=Fortran,style=fortran]{psb_Xasb} routines is specified, data are stored in CPU memory (CSR format).  
To allocate on the GPU, we need to declare mold types with the \lstinline[language=Fortran,style=fortran]{target} attribute:
\begin{lstlisting}[language=Fortran,style=fortran]
type(psb_d_cuda_csrg_sparse_mat), target :: spmold
type(psb_d_vect_cuda),           target :: vmold
type(psb_i_vect_cuda),           target :: imold
\end{lstlisting}
Passing these molds ensures device-resident data and avoids unnecessary host–device transfers.  
The \lstinline[language=Fortran,style=fortran]{mold} attribute in Fortran specifies a prototype or reference object that determines the allocation type and structure of the new variable. When used in conjunction with polymorphic or extensible types, it enables dynamic allocation that conforms to the data representation of the specified mold. Within the state-pattern design adopted by \texttt{PSBLAS}, this mechanism allows the same high-level interface to operate seamlessly across different computational backends (e.g., CPU, CUDA, and different matrix formats). By supplying a mold corresponding to the desired backend, the correct implementation of data structures and associated methods is instantiated at runtime, thereby ensuring both code reusability and backend-specific optimization without altering the user-level algorithmic code.

\paragraph{Preconditioner setup.}
Once the matrix, vector, and communication structures are built, preconditioners can be created:
\begin{lstlisting}[language=Fortran,style=fortran]
type(psb_dprec_type) :: prec
call prec%
\end{lstlisting}
For the AMG methods described in Section~\ref{sec:amg_prec}, the configuration is performed by:
\begin{lstlisting}[language=Fortran,style=fortran,basicstyle=\small]
call prec%
call prec%
call prec%
call prec%
call prec%
call prec%
call prec%
call prec%
call prec%
call prec%
call prec%
\end{lstlisting}
This defines an AMG preconditioner performing one V-cycle, using matching-based aggregation (Sec.~\ref{sec:match_aggreg}) into aggregates of size~8, a degree-8 optimized Chebyshev polynomial smoother~\cite{DAmbra2025} of the first kind, and 30 sweeps of \(\ell_1\)-Jacobi on the coarsest level.  
The preconditioner hierarchy and smoothers are then built via:
\begin{lstlisting}[language=Fortran,style=fortran,basicstyle=\small]
call prec%
call prec%
\end{lstlisting}

\paragraph{Solver invocation.}
After setup, the preconditioner can be used with the available Krylov solvers:
\begin{lstlisting}[language=Fortran,style=fortran,basicstyle=\small]
call psb_krylov('CG', K, prec, f, u, eps, desc, info, &
     itmax=itmax, iter=iter, err=err, itrace=itrace)
\end{lstlisting}
Here, \lstinline{u} is the initial guess on input and the solution on output,  {which readily allows for the use of the all-zero initial guess or, in time-evolving processes, the solution at the previous time step as the initial guess for the current one.} 
\lstinline{eps} sets the tolerance, and \lstinline{itmax} limits iterations.  
The variables \lstinline{iter} and \lstinline{err} return the iteration count and final residual norm, while \lstinline{itrace} controls residual print frequency.

\noindent
This streamlined workflow unifies sparse data distribution, assembly, and solver invocation under a single, portable interface, supporting both CPU and GPU execution.

\section{Numerical experiments}
\label{sec:numerical_experiments}

We consider three test cases based on the 3D Poisson problem~\eqref{eq:poisson}, with discretization performed using the \texttt{GeoPDEs} software~\cite{GeoPDEs1,GeoPDEs2}.\footnote{After installation, these examples can be generated and explored using the \lstinline[language=Matlab]!geopdes_base_examples! routine. We provide brief descriptions here, but refer the reader to the code for full implementation details, including, for instance, the exact placement of control points and weights used to define the NURBS representation in the ring example below.} 
\begin{figure}[htbp]
    \centering
    \includegraphics[width=0.33\columnwidth]{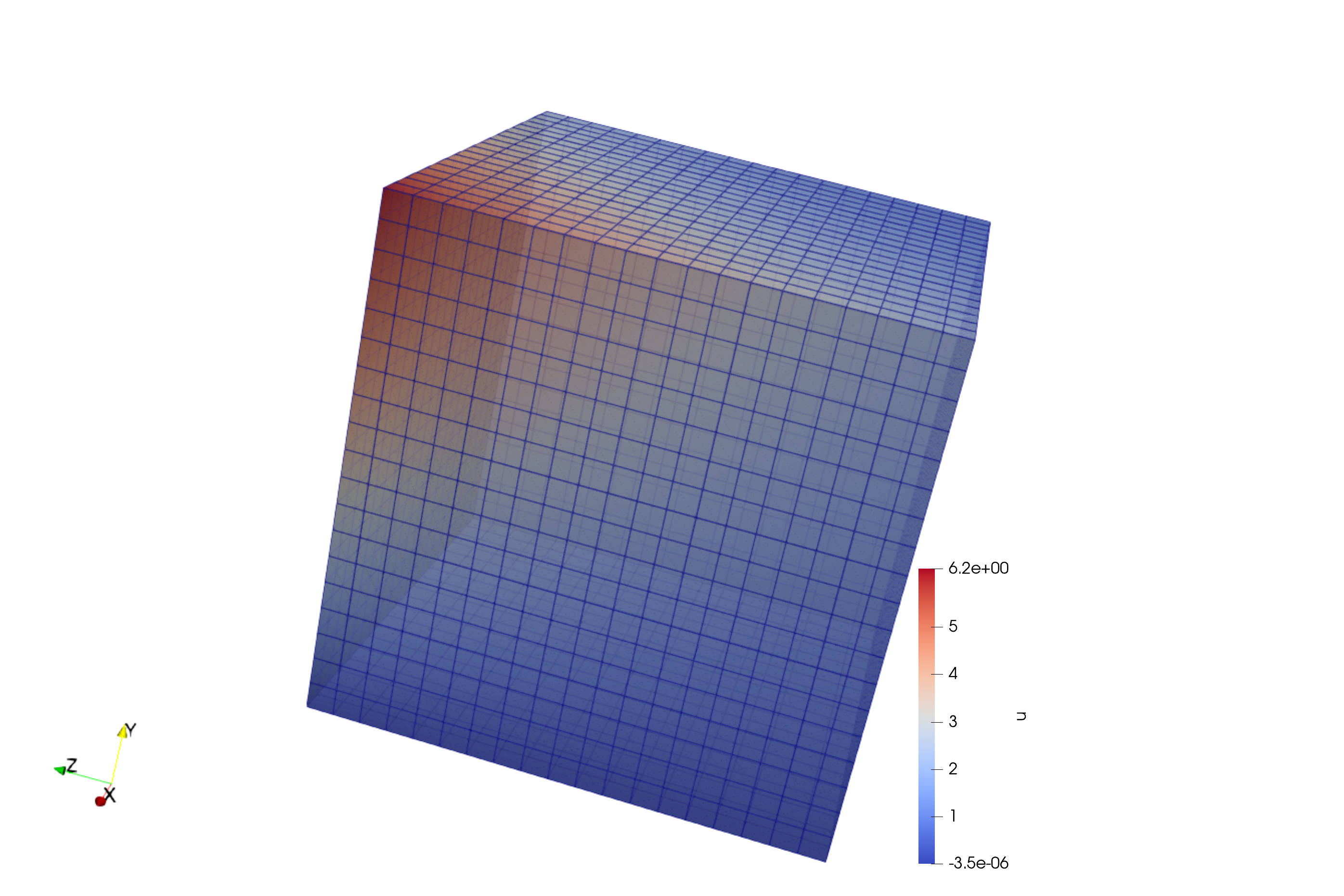}\includegraphics[width=0.33\columnwidth]{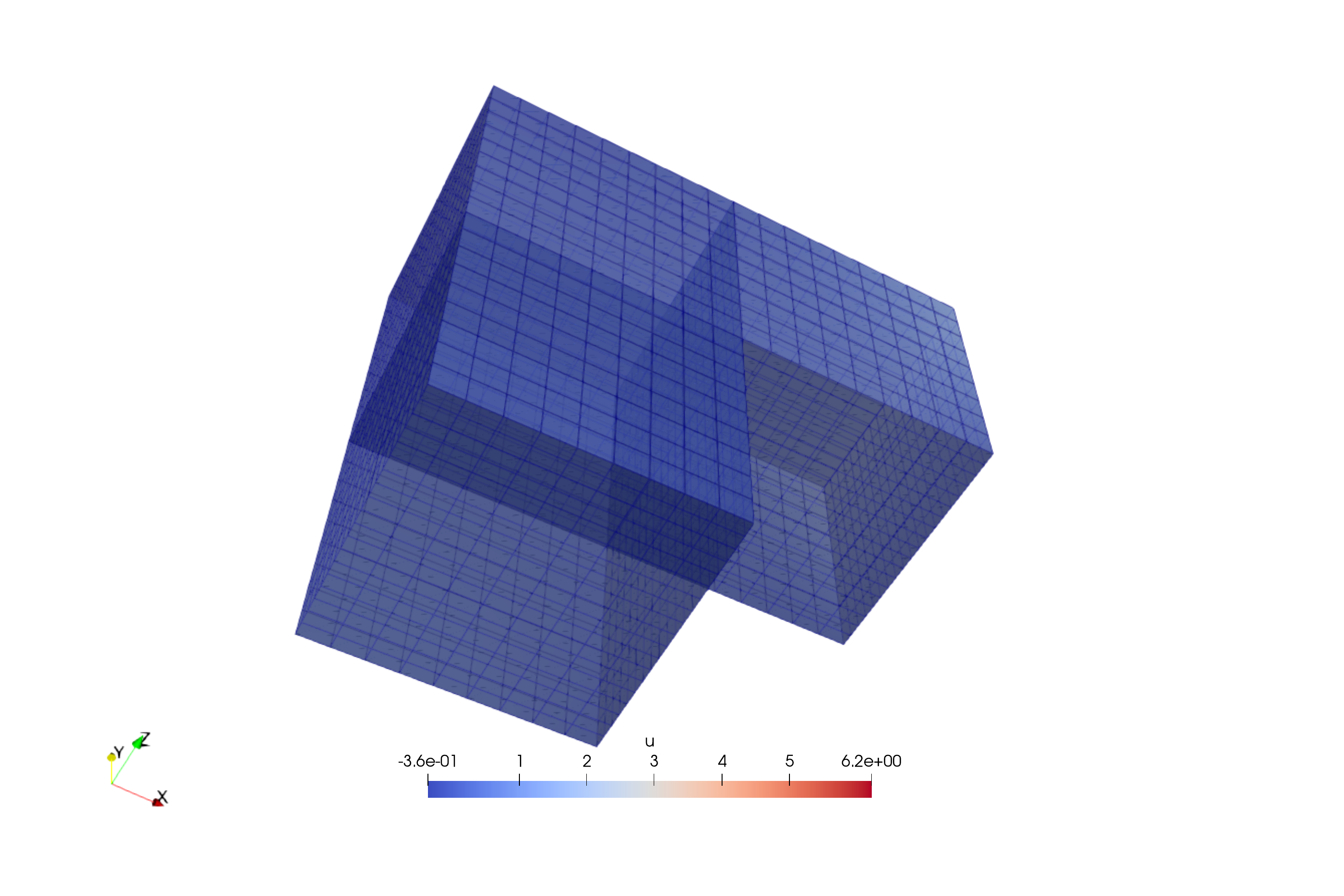}\includegraphics[width=0.33\columnwidth]{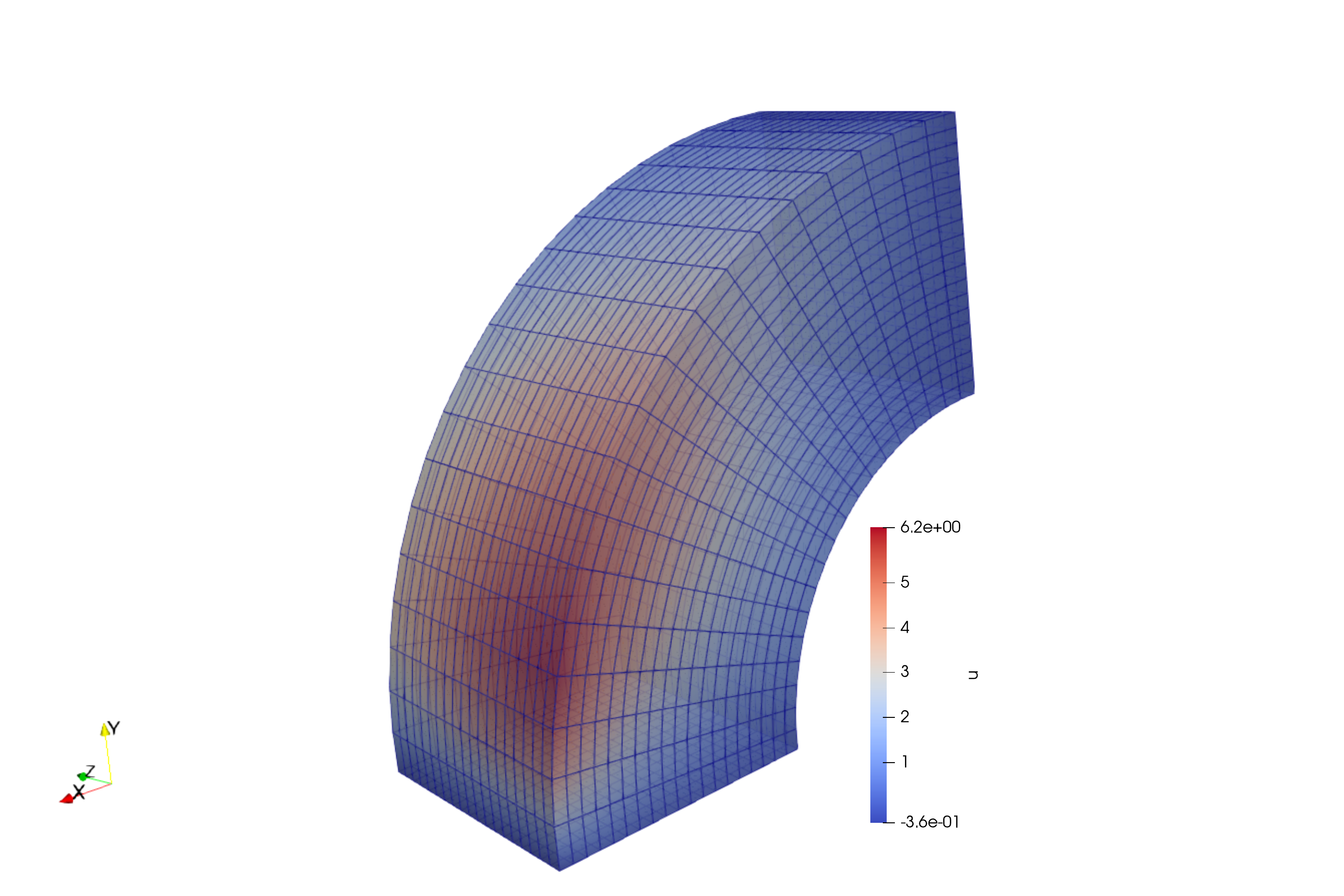}
    \caption{Depiction of the solution for the three test problem considered in the following sections.}
    \label{fig:test_problem}
\end{figure}
The computational domains and exact solutions for these three test cases are shown in Fig.~\ref{fig:test_problem}:
\begin{description}[aaa]
\item[\textbf{Cube}] The computational domain is the unit cube $\Omega = [0,1]^3$. Dirichlet boundary conditions are imposed on three of its faces, while Neumann boundary conditions are applied on the remaining three, as described below:
\[
\begin{split}
f(x,y,z) &= -\exp(x+z)\sin(y), \qquad g_D(x,y,z) = \exp(x+z)\sin(y),\\
g_N(x,y,z) &= 
\begin{cases}
\exp(x+z) \cos(y) & \text{if } \text{side} = 4,\\
- \exp(x+z) \sin(y) & \text{if } \text{side} = 5,\\
\exp(x+z) \sin(y) & \text{if } \text{side} = 6.
\end{cases}
\end{split}
\]

\item[\textbf{L-shaped}] The domain is a three-dimensional L-shaped region constructed from three patches (see Section~\ref{sec:multipatch}), defined as the union $\Omega = [0,1]^3 \cup [0,-1]^3 \cup [1,2]^3$. Dirichlet boundary conditions are imposed on six faces of the domain, while Neumann conditions are applied to the remaining two. The source term and boundary data are given by:
\begin{equation} \label{eq:Lshapedcoeff}
\begin{split}
f(x,y,z) &= \exp(x)\cos(z)\left[-2\cos(xy)y + \sin(xy)(y^2 + x^2)\right], \\
g_D(x,y,z) &= \exp(x)\sin(xy)\cos(z), \\
g_N(x,y,z) &=
\begin{cases}
- x \exp(x) \cos(xy) \cos(z) & \text{if } \text{side} = 3,\\
\cos(z) \exp(x) \left[\sin(xy) + y \cos(xy)\right] & \text{if } \text{side} = 4,\\
- \cos(z) \exp(x) \left[\sin(xy) + y \cos(xy)\right] & \text{if } \text{side} = 5,\\
x \exp(x) \cos(xy) \cos(z) & \text{if } \text{side} = 6,\\
\exp(x) \sin(xy) \sin(z) & \text{if } \text{side} = 7,\\
- \exp(x) \sin(xy) \sin(z) & \text{if } \text{side} = 8.
\end{cases}
\end{split}
\end{equation}

\item[\boldmath$\nicefrac{1}{4}$\textbf{-ring}] The domain is a thick quarter-ring centered at the origin, with inner radius 1, outer radius 2, and height 1, entirely contained within the positive orthant. The geometry is modeled using a NURBS representation, while the equation is discretized using B-spline basis functions, resulting in a non-isoparametric approach (see Section~\ref{sec:iso-or-noniso}). Dirichlet conditions are imposed on the first three boundary faces, and Neumann conditions on the remaining three. The functions $f$ and $g_D$ are as defined in~\eqref{eq:Lshapedcoeff}, while the Neumann data is given by:
\begingroup
\allowdisplaybreaks
\[
g_N(x, y, z) = 
\begin{cases}
- \exp(x) \cos(z) \left[\sin(xy) + y \cos(xy)\right] & \text{if } \text{side} = 4, \\
\exp(x) \sin(xy) \sin(z) & \text{if } \text{side} = 5, \\
- \exp(x) \sin(xy) \sin(z) & \text{if } \text{side} = 6.
\end{cases}
\]
\endgroup
\end{description}

In all test cases, we adopt maximum-regularity B-spline bases and export the resulting system matrices and right-hand side vectors in Matrix Market format.

Then \texttt{PSCToolkit} is used to read the data from file, distribute them among the computational processes via the \texttt{METIS} library~\cite{METIS} and solve the resulting linear system by means of the flexible variant of the CG algorithm available through PSBLAS with a request on the relative tolerance of $10^{-6}$, {the initial guess is always chosen as the $\mathbf{0}$ vector}. 

Numerical experiments have been run on the nodes of the Amelia cluster at the Institute for Applied Computing ``Mauro Picone'' of the National Research Council of Italy in Naples. Each node is equipped with a Intel\textsuperscript{\textregistered} Xeon\textsuperscript{\textregistered} Gold 6338 CPU at 2.00~\si{\giga\Hz} (1 thread per core, 32 cores per socket and 2 sockets), 1~\si{\tera\byte} of RAM, and 4 NVIDIA A30 GPUs with 24 \si{\giga\byte} of HBM2 memory each. The software environment consists of \texttt{gcc} v12.2.1, CUDA toolkit v12.9 (with \texttt{nvcc} compiler v12.9.41), \texttt{openmpi} v4.1.6, and \texttt{openblas} v0.3.29 (Skylake-X optimized), together with the numerical libraries \texttt{PSBLAS}~v3.9.0-rc3 and \texttt{AMG4PSBLAS}~v1.2-rc3. In all the following examples, we restrict the number of MPI tasks per node to 32 for pure MPI executions, and to 4 when enabling GPU acceleration. In the latter case, each MPI task is associated with a distinct GPU device on the node, ensuring an efficient one-to-one mapping between MPI ranks and GPUs.

\subsection{Cube}\label{sec:cube}
We begin by considering the Poisson problem on the unit cube. Our initial focus is on evaluating the algorithmic scalability of the matching-based AMG preconditioner (V-Cycle, smoothed, aggregates of size 8) with 1\textsuperscript{st}-kind optimized Chebyshev accelerated smoother ($\ell_1$-Jacobi) and the CG preconditioned by a single sweep of weighted Jacobi as coarse solver set to achieve a tolerance of $10^{-4}$ or stop in $30$ iterations. {We note that the aggregate size is chosen to balance the regularity properties of the smoothed prolongator and the coarsening ratio, thereby keeping the overall hierarchy complexity sufficiently low.} The number of iterations of the smoothers, i.e., the degree of the Chebyshev polynomial
$\deg_p$ and, equivalently, the number of matrix–vector products, is selected according to
the spline degree $p$ of the basis. In particular, we take
$\deg_3 = 8$, $\deg_4 = 12$, $\deg_5 = 14$, and $\deg_6 = 16$.
This choice is motivated by the spectral behavior of stiffness matrices described in
\cite[Lemma~4.1]{Manni1}: as the spline degree $p$ increases, the normalized spectral distribution develops an increasingly sharp quadratic zero at $\theta = 0$ and decays only slowly toward 
its mirror point $\theta = \pi$, where its value decreases exponentially with $p$.
Consequently, higher spline degrees generate more persistent low–frequency components that are
poorly damped by simple relaxation schemes. To achieve $p$–robust smoothing efficiency,
heavier (higher–degree) Chebyshev smoothers are therefore required.
 We analyze the number of iterations required for convergence as a function of the number of subdivisions $k$ along each edge of the cube and the polynomial degree $p$ of the B-spline basis functions, while enforcing the maximum possible continuity, namely $\mathcal{C}^{p-1}$.
\begin{table}[htbp]
    \centering
    \begin{subtable}{0.2\columnwidth}
    \begin{tabular}{l|cccc}
    \toprule
    $k \backslash p$  &  3 & 4 & 5 & 6\\
    \midrule
12              & 10              & 22              & 27              & 25             \\
24              & 11              & 23              & 29              & 32             \\
48              & 11              & 22              & 27              & 34             \\
96              & 11              & 21              & $\dagger$       & $\dagger$      \\
    \bottomrule
    \end{tabular}
    \caption{Iteration count}
    \end{subtable}\hfil
    \begin{subtable}{0.4\columnwidth}
    \begin{tabular}{l|cccc}
    \toprule
    $k \backslash p$  &  3 & 4 & 5 & 6\\
    \midrule
12              & 2730            & 3360            & 4080            & 4896            \\
24              & 17550           & 19656           & 21924           & 24360           \\
48              & 124950          & 132600          & 140556          & 148824          \\
96              & 941094          & 970200          & $\dagger$       & $\dagger$      \\
    \bottomrule
    \end{tabular}
    \caption{Matrix size}
    \end{subtable}\hfil
    \begin{subtable}{0.25\columnwidth}
    \begin{tabular}{l|cccc}
    \toprule
    $k \backslash p$  &  3 & 4 & 5 & 6\\
    \midrule
12              & 1.20 & 1.12 & 1.09 & 1.07 \\
24              & 1.23 & 1.19 & 1.15 & 1.12 \\
48              & 1.29 & 1.25 & 1.21 & 1.17 \\
96              & 1.34 & 1.30 & $\dagger$       & $\dagger$      \\
    \bottomrule
    \end{tabular}
    \caption{Operator complexity}
    \end{subtable}
    
    \caption{Result for the Poisson on the cube case on a single task for varying number of subdivision $k$ of the base mesh, and degree $p$ of the B-Spline basis. The result highlighted with a $\dagger$ are due to the sequential matrix assembler not having enough RAM to build the discretization on a single node. }
    \label{tab:poisson_cube}
\end{table}
To assess parallel efficiency, we conduct a strong scaling experiment. In this configuration, the problem size is kept fixed while the number of MPI tasks is increased from $1$ to $512$ employing at most $32$ tasks per node. 

We begin by examining the behavior with respect to operator complexity. As discussed in Section~\ref{sec:match_aggreg}, the aggregation algorithm relies on an approximate graph matching. However, since the graph is now distributed across multiple tasks, this matching differs from the non-distributed case~\cite{6009071}.
\begin{figure}[htbp]
    \centering
    \includegraphics[width=\columnwidth]{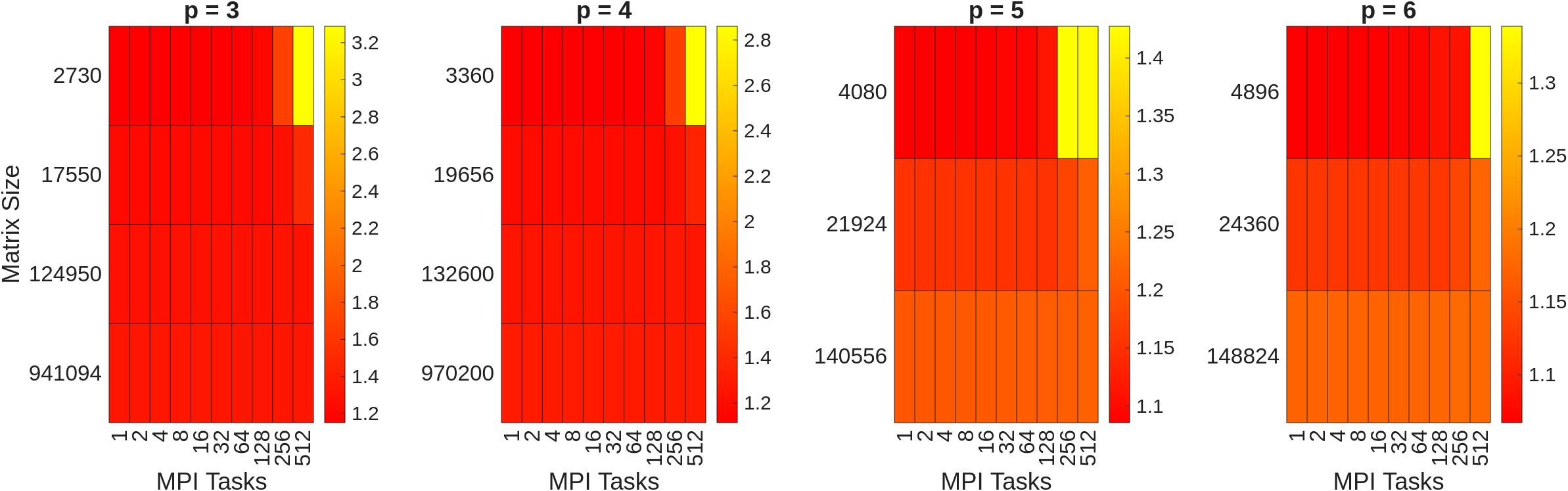}
    \caption{Operator complexity for the Poisson cube case for increasing values of $k$ and $p$, and number of MPI tasks running from $1$ to $512$.}
    \label{fig:cube_opc}
\end{figure}
The results closely match the single--process behaviour reported in
Table~\ref{tab:poisson_cube}.  
For small local problem sizes, however, each MPI task owns only a very
limited portion of the graph; as a consequence, the effective coarsening
ratio deteriorates, reducing the ability of the aggregation procedure to
form aggregates of the target size.  
This degradation is further amplified by the coarsening stopping
criterion, which in our implementation enforces a minimum global coarse
size proportional to the number of tasks (specifically, the coarsening
phase is terminated once the global coarse size drops below $50*ntasks$).

\begin{figure}[hbt]
    \definecolor{mycolor1}{rgb}{0.06667,0.44314,0.74510}%
\definecolor{mycolor2}{rgb}{0.86667,0.32941,0.00000}%
\definecolor{mycolor3}{rgb}{0.12941,0.12941,0.12941}%
\begin{tikzpicture}

\begin{axis}[%
width=0.36\columnwidth,
height=0.196\columnwidth,
at={(0\columnwidth,0.27\columnwidth)},
scale only axis,
xmode=log,
xmin=1,
xmax=512,
xtick={1,2,4,8,16,32,64,128,256,512},
xticklabels={1,2,4,8,16,32,64,128,256,512},
xticklabel style={rotate=45,font=\footnotesize,anchor=east},
xminorticks=true,
xlabel style={font=\color{mycolor3}},
separate axis lines,
every outer y axis line/.append style={mycolor1},
every y tick label/.append style={font=\color{mycolor1}},
every y tick/.append style={mycolor1},
ymode=log,
ymin=0.5,
ymax=200,
yminorticks=true,
ylabel style={font=\color{mycolor1}},
ylabel={Time (s)},
axis background/.style={fill=white},
title style={font=\small},
title={$p = 3$},
title style={yshift=-1em},
yticklabel pos=left,
xmajorgrids,
xminorgrids,
ymajorgrids,
yminorgrids
]
\addplot [color=mycolor1, line width=2.0pt, mark=o, mark options={solid, mycolor1}, forget plot]
  table[row sep=crcr]{%
512	1.15216\\
256	1.19071\\
128	1.64733\\
64	2.8398\\
32	5.36222\\
16	9.21073\\
8	18.6799\\
4	34.1878\\
2	67.3164\\
1	132.37\\
};
\end{axis}
\begin{axis}[%
width=0.36\columnwidth,
height=0.196\columnwidth,
at={(0\columnwidth,0.27\columnwidth)},
scale only axis,
xmode=log,
xmin=1,
xmax=512,
xtick={1,2,4,8,16,32,64,128,256,512},
xticklabels={},
xminorticks=true,
xlabel style={font=\color{mycolor3}},
separate axis lines,
every outer y axis line/.append style={mycolor2},
every y tick label/.append style={font=\color{mycolor2}},
every y tick/.append style={mycolor2},
ymin=10,
ymax=12,
yminorticks=true,
ylabel
axis background/.style={fill=none},
title style={font=\bfseries\color{mycolor3}},
yticklabel pos=right,
]
\addplot [color=mycolor2, dashed, line width=2.0pt, forget plot]
  table[row sep=crcr]{%
512	11\\
256	11\\
128	11\\
64	11\\
32	11\\
16	11\\
8	11\\
4	11\\
2	11\\
1	11\\
};
\end{axis}

\begin{axis}[%
width=0.36\columnwidth,
height=0.196\columnwidth,
at={(0.44\columnwidth,0.27\columnwidth)},
scale only axis,
xmode=log,
xmin=1,
xmax=512,
xtick={1,2,4,8,16,32,64,128,256,512},
xticklabels={1,2,4,8,16,32,64,128,256,512},
xticklabel style={rotate=45,font=\small,anchor=east},
xminorticks=true,
xlabel style={font=\color{mycolor3}},
separate axis lines,
every outer y axis line/.append style={mycolor1},
every y tick label/.append style={font=\color{mycolor1}},
every y tick/.append style={mycolor1},
ymin=2,
ymax=700,
ymode=log,
yminorticks=true,
ylabel style={font=\color{mycolor1}},
axis background/.style={fill=white},
title style={font=\small},
title={$p = 4$},
title style={yshift=-1em},
yticklabel pos=left,
xmajorgrids,
xminorgrids,
ymajorgrids,
yminorgrids
]
\addplot [color=mycolor1, line width=2.0pt, mark=o, mark options={solid, mycolor1}, forget plot]
  table[row sep=crcr]{%
512	3.91637\\
256	4.87987\\
128	7.82285\\
64	14.8386\\
32	26.8899\\
16	44.9079\\
8	82.8711\\
4	162.286\\
2	316.376\\
1	618.677\\
};
\end{axis}
\begin{axis}[%
width=0.36\columnwidth,
height=0.196\columnwidth,
at={(0.44\columnwidth,0.27\columnwidth)},
scale only axis,
xmode=log,
xmin=1,
xmax=512,
xtick={1,2,4,8,16,32,64,128,256,512},
xticklabels={},
xminorticks=true,
separate axis lines,
every outer y axis line/.append style={mycolor2},
every y tick label/.append style={font=\color{mycolor2}},
every y tick/.append style={mycolor2},
ymin=16,
ymax=20,
yminorticks=true,
axis background/.style={fill=none},
yticklabel pos=right,
ylabel style={font=\color{mycolor2}},
ylabel={Iterations},
]
\addplot [color=mycolor2, dashed, line width=2.0pt, forget plot]
  table[row sep=crcr]{%
512	18\\
256	18\\
128	19\\
64	19\\
32	18\\
16	17\\
8	17\\
4	17\\
2	17\\
1	17\\
};
\end{axis}

\begin{axis}[%
width=0.36\columnwidth,
height=0.196\columnwidth,
at={(0\columnwidth,0\columnwidth)},
scale only axis,
xmode=log,
xmin=1,
xmax=512,
xtick={1,2,4,8,16,32,64,128,256,512},
xticklabels={1,2,4,8,16,32,64,128,256,512},
xticklabel style={rotate=45,font=\small,anchor=east},
xminorticks=true,
xlabel style={font=\color{mycolor3}},
xlabel={MPI Tasks},
separate axis lines,
every outer y axis line/.append style={mycolor1},
every y tick label/.append style={font=\color{mycolor1}},
every y tick/.append style={mycolor1},
ymode=log,
ymin=2,
ymax=300,
yminorticks=true,
ylabel style={font=\color{mycolor1}},
ylabel={Time (s)},
axis background/.style={fill=white},
title style={font=\small},
title={$p = 5$},
title style={yshift=-1em},
yticklabel pos=left,
xmajorgrids,
xminorgrids,
ymajorgrids,
yminorgrids
]
\addplot [color=mycolor1, line width=2.0pt, mark=o, mark options={solid, mycolor1}, forget plot]
  table[row sep=crcr]{%
512	3.27363\\
256	3.40639\\
128	3.88562\\
64	6.30663\\
32	10.7437\\
16	18.7867\\
8	35.1372\\
4	66.2161\\
2	122.294\\
1	242.847\\
};
\end{axis}
\begin{axis}[%
width=0.36\columnwidth,
height=0.196\columnwidth,
at={(0\columnwidth,0\columnwidth)},
scale only axis,
xmode=log,
xmin=1,
xmax=512,
xtick={1,2,4,8,16,32,64,128,256,512},
xticklabels={},
xminorticks=true,
separate axis lines,
every outer y axis line/.append style={mycolor2},
every y tick label/.append style={font=\color{mycolor2}},
every y tick/.append style={mycolor2},
ymin=22,
ymax=28,
yminorticks=true,
axis background/.style={fill=none},
yticklabel pos=right,
]
\addplot [color=mycolor2, dashed, line width=2.0pt, forget plot]
  table[row sep=crcr]{%
512	23\\
256	26\\
128	26\\
64	26\\
32	25\\
16	26\\
8	26\\
4	26\\
2	25\\
1	25\\
};
\end{axis}

\begin{axis}[%
width=0.36\columnwidth,
height=0.196\columnwidth,
at={(0.44\columnwidth,0\columnwidth)},
scale only axis,
xmode=log,
xmin=1,
xmax=512,
xtick={  1,   2,   4,   8,  16,  32,  64, 128, 256, 512},
xticklabels={  1,   2,   4,   8,  16,  32,  64, 128, 256, 512},
x tick label style={rotate=45,font=\small},
xminorticks=true,
xlabel style={font=\color{mycolor3}},
xlabel={MPI Tasks},
separate axis lines,
every outer y axis line/.append style={mycolor1},
every y tick label/.append style={font=\color{mycolor1}\small},
every y tick/.append style={mycolor1},
ymin=5,
ymax=700,
ymode=log,
axis background/.style={fill=none},
title style={font=\small\color{mycolor3}},
title={$p = 6$},
title style={yshift=-1em},
yticklabel pos=left,
xmajorgrids,
xminorgrids,
ymajorgrids,
yminorgrids
]
\addplot [color=mycolor1, line width=2.0pt, mark=o, mark options={solid, mycolor1}, forget plot]
  table[row sep=crcr]{%
512	8.49992\\
256	7.81247\\
128	9.75418\\
64	15.7963\\
32	28.9733\\
16	48.6326\\
8	86.5748\\
4	163.819\\
2	315.153\\
1	620.898\\
};
\end{axis}
\begin{axis}[%
width=0.36\columnwidth,
height=0.196\columnwidth,
at={(0.44\columnwidth,0\columnwidth)},
scale only axis,
xmode=log,
xmin=1,
xmax=512,
xtick={  1,   2,   4,   8,  16,  32,  64, 128, 256, 512},
xticklabels={  1,   2,   4,   8,  16,  32,  64, 128, 256, 512},
x tick label style={rotate=45,font=\small},
xminorticks=true,
xlabel style={font=\color{mycolor3}},
xlabel={MPI Tasks},
separate axis lines,
every outer y axis line/.append style={mycolor2},
every y tick label/.append style={font=\color{mycolor2}\small},
every y tick/.append style={mycolor2},
ymin=33,
ymax=36,
axis background/.style={fill=none},
ylabel style={font=\color{mycolor2}},
ylabel={Iterations},
yticklabel pos=right,
xmajorgrids,
xminorgrids,
ymajorgrids,
yminorgrids
]
\addplot [color=mycolor2, dashed, line width=2.0pt, forget plot]
  table[row sep=crcr]{%
512	35\\
256	35\\
128	35\\
64	35\\
32	35\\
16	34\\
8	34\\
4	34\\
2	34\\
1	34\\
};
\end{axis}

\end{tikzpicture}%
    \caption{Strong scaling for the Poisson on a 3D Cube problem for a number of task running from $1$ to $512$ {on the largest available test case for each $p$, i.e., sizes $941094$, $970200$,  $140556$ and $148824$ respectively.}}
    \label{fig:cube_h_strong_scaling}
\end{figure}
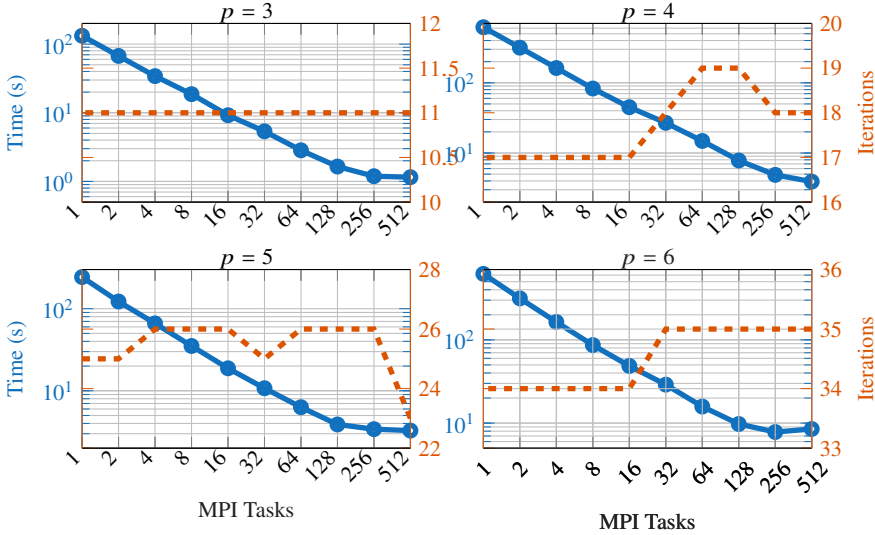

Figure~\ref{fig:cube_h_strong_scaling} reports the wall-clock solution time (left column) and the corresponding iteration count (right column) as the number of MPI tasks increases, for polynomial degrees $p=3,4,5,6$. Across all polynomial degrees, we observe a significant reduction in time-to-solution, demonstrating good strong scaling up to $512$ tasks. Throughout most of the scaling range, the iteration count remains nearly constant, confirming that the preconditioner maintains its algorithmic robustness under parallel execution.
Examining the strong-scaling performance quantitatively, for $p = 3$ the solve time decreases from $132.37\,\mathrm{s}$ on a single task to $1.15\,\mathrm{s}$ on 512 tasks, yielding a speedup of approximately $115\times$ and a parallel efficiency of $22\%$. For $p = 4$, the time drops from $618.68\,\mathrm{s}$ to $3.92\,\mathrm{s}$, corresponding to a speedup of $158\times$ and a parallel efficiency of $31\%$. The $p = 5$ case shows comparable scalability, with the time decreasing from $242.85\,\mathrm{s}$ to $3.27\,\mathrm{s}$, resulting in a speedup of $74\times$ and a parallel efficiency of $14\%$. For $p = 6$, the time reduces from $620.90\,\mathrm{s}$ to $8.50\,\mathrm{s}$, giving a speedup of $73\times$ and a parallel efficiency of $14\%$. These results indicate that the preconditioner scales well in a distributed-memory environment, although parallel efficiencies at 512 tasks remain moderate (14--31\%), reflecting the increasing impact of communication overhead for smaller local problem sizes. At the highest processor counts, particularly for larger values of $p$, a slight increase in the iteration count is observed; however, this increase is modest, and the overall time-to-solution continues to decrease, indicating that performance remains dominated by computational workload rather than convergence degradation.

\subsection{Multipatch \texorpdfstring{$L$-shaped domain}{L-shaped domain}}\label{sec:multipatch-experiments}
We now turn to the three-patch $L$-shaped domain (cf. the problem definition in Section~\ref{sec:numerical_experiments}). %
Table~\ref{tab:poisson_lshape} summarizes the single-process behavior: iteration counts, matrix sizes, and operator complexity for increasing mesh subdivisions $k$ and B-spline degree $p$. The counts remain modest and only mildly sensitive to $k$, confirming that the preconditioner maintains good algorithmic scalability also in the presence of patch interfaces. As expected, the operator complexity is slightly higher than in the single-patch cube due to the additional interface DoFs and the reduced inter-patch continuity, which increase the aggregate connectivity.
\begin{table}[htbp]
    \centering
    \begin{subtable}{0.2\columnwidth}
    \begin{tabular}{l|ccccc}
    \toprule
    $k \backslash p$  &  2 & 3 & 4 & 5 & 6 \\
    \midrule
12              & 6               & 8               & 15             & 20             & 23            \\
24              & 7               & 9               & 14              & 17              & 27             \\
48              & 7               & 9               & 13              & 15              & 16             \\
96              & 10              & 10              & $\dagger$              & $\dagger$              & $\dagger$              \\
    \bottomrule
    \end{tabular}
    \caption{Iteration count}
      \end{subtable}\hfil
    \begin{subtable}{0.42\columnwidth}
    \begin{tabular}{l|ccccc}
    \toprule
    $k \backslash p$  &  2 & 3 & 4 & 5 & 6 \\
    \midrule
12              & 5772            & 7280            & 9030             & 11024             & 13328             \\
24              & 43800           & 49400           & 55485            & 61922           & 69020           \\
48              & 341040          & 362600          & 385050           & 408408           & 432692           \\
96              & 2691168         & 2775752         & $\dagger$          & $\dagger$          & $\dagger$          \\
    \bottomrule
    \end{tabular}
    \caption{Matrix size}
    \end{subtable}
    
    \begin{subtable}{0.25\columnwidth}
    \begin{tabular}{l|ccccc}
    \toprule
    $k \backslash p$  &  2 & 3 & 4 & 5 & 6 \\
    \midrule
12              & 1.35 & 1.20 & 1.16 & 1.13 & 1.11 \\
24              & 1.49 & 1.27 & 1.23 & 1.19 & 1.16 \\
48              & 1.59 & 1.32 & 1.28 & 1.23 & 1.19 \\
96              & 1.65 & 1.35 & $\dagger$ & $\dagger$ & $\dagger$ \\
    \bottomrule
    \end{tabular}
    \caption{Operator complexity}
    \end{subtable}
  
        \caption{Result for the Poisson on the multipatch $L$-shaped domain case on a single task for varying number of subdivision $k$ of the base mesh, and degree $p$ of the B-Spline basis. The result highlighted with a $\dagger$ are due to the sequential matrix assembler not having enough RAM to build the discretization on a single node. }
    \label{tab:poisson_lshape}
\end{table}
\begin{figure}[htbp]
    \centering
    \includegraphics[width=\columnwidth]{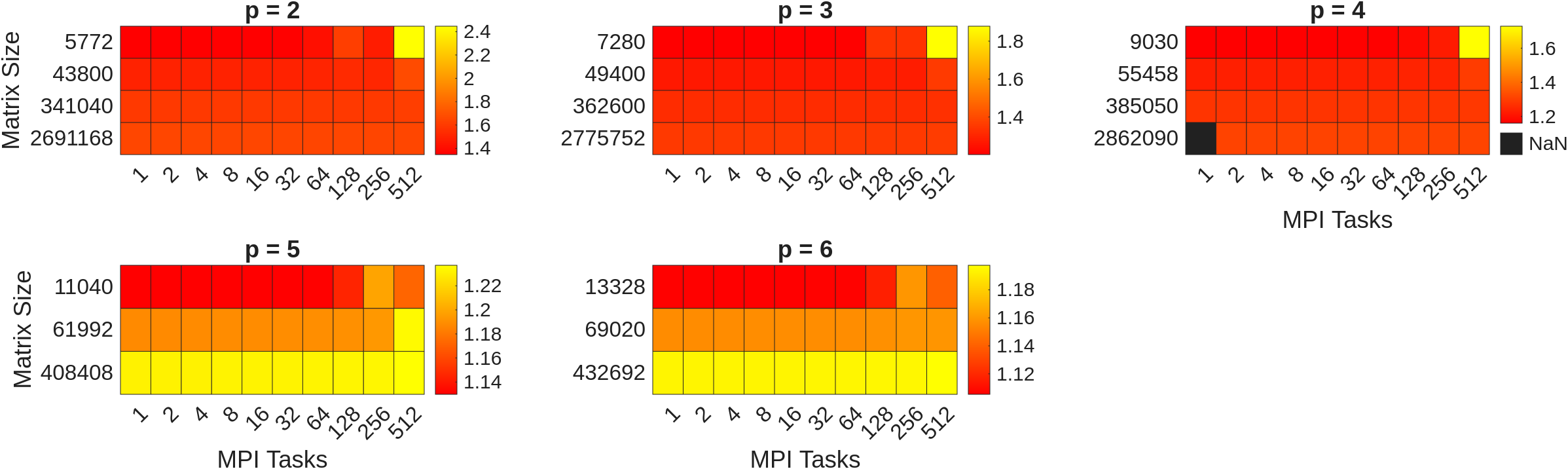}
    \caption{Operator complexity for the multipatch $L$-shaped case for increasing values of $k$ and $p$, and number of MPI tasks running from $1$ to $512$.}
    \label{fig:multipatch_opc}
\end{figure}
We next investigate parallel efficiency. As in the Cube study, we perform a strong-scaling experiment in which the problem size is fixed while the number of MPI tasks increases from $1$ to $512$ (up to $32$ tasks per node). Figure~\ref{fig:multipatch_opc} reports the variation of operator complexity with the task count: for small matrices, the complexity grows as the number of tasks becomes large. As in the cube case, this behavior arises because very small local partitions prevent the formation of aggregates of size~8, leading to smaller aggregates (and thus more coarse unknowns) and increasing the operator complexity. Patch interfaces introduce additional connectivity and slightly amplify this effect, but the increase is primarily driven by the aggregate-size constraint rather than by reduced matching quality.

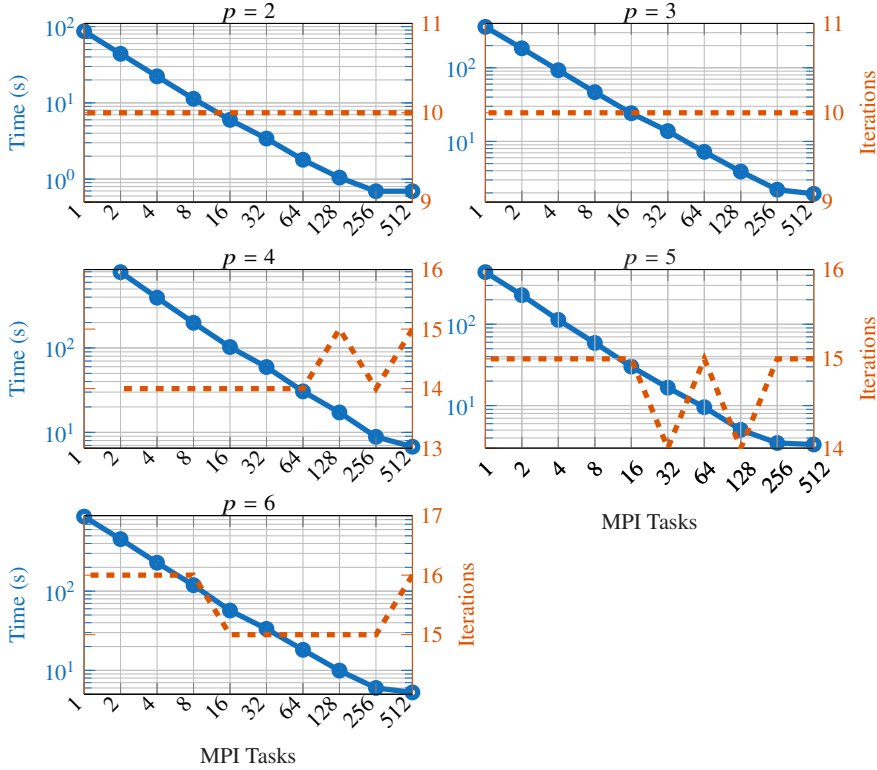
\begin{figure}
    \sidecaption
    \definecolor{mycolor1}{rgb}{0.06667,0.44314,0.74510}%
\definecolor{mycolor2}{rgb}{0.86667,0.32941,0.00000}%
\definecolor{mycolor3}{rgb}{0.12941,0.12941,0.12941}%
\begin{tikzpicture}

\begin{axis}[%
width=0.36\columnwidth,
height=0.196\columnwidth,
at={(0\columnwidth,0.27\columnwidth)},
scale only axis,
xmode=log,
xmin=1,
xmax=512,
xtick={1,2,4,8,16,32,64,128,256,512},
xticklabels={1,2,4,8,16,32,64,128,256,512},
xticklabel style={rotate=45,font=\footnotesize,anchor=east},
xminorticks=true,
xlabel style={font=\color{mycolor3}},
separate axis lines,
every outer y axis line/.append style={mycolor1},
every y tick label/.append style={font=\color{mycolor1}},
every y tick/.append style={mycolor1},
ymode=log,
ymin=0.5,
ymax=110,
yminorticks=true,
ylabel style={font=\color{mycolor1}},
ylabel={Time (s)},
axis background/.style={fill=white},
title style={font=\small},
title={$p = 2$},
title style={yshift=-1em},
yticklabel pos=left,
xmajorgrids,
xminorgrids,
ymajorgrids,
yminorgrids
]
\addplot [color=mycolor1, line width=2.0pt, mark=o, mark options={solid, mycolor1}, forget plot]
  table[row sep=crcr]{%
512	0.693995\\
256	0.693945\\
128	1.04701\\
64	1.79453\\
32	3.40724\\
16	5.98823\\
8	11.325\\
4	22.2544\\
2	43.9917\\
1	87.091\\
};
\end{axis}
\begin{axis}[%
width=0.36\columnwidth,
height=0.196\columnwidth,
at={(0\columnwidth,0.27\columnwidth)},
scale only axis,
xmode=log,
xmin=1,
xmax=512,
xtick={1,2,4,8,16,32,64,128,256,512},
xticklabels={},
xminorticks=true,
xlabel style={font=\color{mycolor3}},
separate axis lines,
every outer y axis line/.append style={mycolor2},
every y tick label/.append style={font=\color{mycolor2}},
every y tick/.append style={mycolor2},
ymin=9,
ymax=11,
yminorticks=true,
ytick={9,10,11},
axis background/.style={fill=none},
title style={font=\bfseries\color{mycolor3}},
yticklabel pos=right,
]
\addplot [color=mycolor2, dashed, line width=2.0pt, forget plot]
  table[row sep=crcr]{%
512	10\\
256	10\\
128	10\\
64	10\\
32	10\\
16	10\\
8	10\\
4	10\\
2	10\\
1	10\\
};
\end{axis}

\begin{axis}[%
width=0.36\columnwidth,
height=0.196\columnwidth,
at={(0.44\columnwidth,0.27\columnwidth)},
scale only axis,
xmode=log,
xmin=1,
xmax=512,
xtick={1,2,4,8,16,32,64,128,256,512},
xticklabels={1,2,4,8,16,32,64,128,256,512},
xticklabel style={rotate=45,font=\small,anchor=east},
xminorticks=true,
xlabel style={font=\color{mycolor3}},
separate axis lines,
every outer y axis line/.append style={mycolor1},
every y tick label/.append style={font=\color{mycolor1}},
every y tick/.append style={mycolor1},
ymode=log,
ymin=1.5,
ymax=400,
ymode=log,
yminorticks=true,
ylabel style={font=\color{mycolor1}},
axis background/.style={fill=white},
title style={font=\small},
title={$p = 3$},
title style={yshift=-1em},
yticklabel pos=left,
xmajorgrids,
xminorgrids,
ymajorgrids,
yminorgrids
]
\addplot [color=mycolor1, line width=2.0pt, mark=o, mark options={solid, mycolor1}, forget plot]
  table[row sep=crcr]{%
512	1.94327\\
256	2.21041\\
128	3.89778\\
64	7.21018\\
32	13.818\\
16	24.148\\
8	46.6713\\
4	92.6605\\
2	183.949\\
1	359.347\\
};
\end{axis}
\begin{axis}[%
width=0.36\columnwidth,
height=0.196\columnwidth,
at={(0.44\columnwidth,0.27\columnwidth)},
scale only axis,
xmode=log,
xmin=1,
xmax=512,
xtick={1,2,4,8,16,32,64,128,256,512},
xticklabels={},
xminorticks=true,
separate axis lines,
every outer y axis line/.append style={mycolor2},
every y tick label/.append style={font=\color{mycolor2}},
every y tick/.append style={mycolor2},
ymin=9,
ymax=11,
yminorticks=true,
ytick={9,10,11},
axis background/.style={fill=none},
yticklabel pos=right,
ylabel style={font=\color{mycolor2}},
ylabel={Iterations},
]
\addplot [color=mycolor2, dashed, line width=2.0pt, forget plot]
  table[row sep=crcr]{%
512	10\\
256	10\\
128	10\\
64	10\\
32	10\\
16	10\\
8	10\\
4	10\\
2	10\\
1	10\\
};
\end{axis}

\begin{axis}[%
width=0.36\columnwidth,
height=0.196\columnwidth,
at={(0\columnwidth,0\columnwidth)},
scale only axis,
xmode=log,
xmin=1,
xmax=512,
xtick={1,2,4,8,16,32,64,128,256,512},
xticklabels={1,2,4,8,16,32,64,128,256,512},
xticklabel style={rotate=45,font=\small,anchor=east},
xminorticks=true,
xlabel style={font=\color{mycolor3}},
separate axis lines,
every outer y axis line/.append style={mycolor1},
every y tick label/.append style={font=\color{mycolor1}},
every y tick/.append style={mycolor1},
ymode=log,
ymin=6.5,
ymax=850,
yminorticks=true,
ylabel style={font=\color{mycolor1}},
ylabel={Time (s)},
axis background/.style={fill=white},
title style={font=\small},
title={$p = 4$},
title style={yshift=-1em},
yticklabel pos=left,
xmajorgrids,
xminorgrids,
ymajorgrids,
yminorgrids
]
\addplot [color=mycolor1, line width=2.0pt, mark=o, mark options={solid, mycolor1}, forget plot]
  table[row sep=crcr]{%
512	6.75947\\
256	8.85094\\
128	17.2197\\
64	30.5903\\
32	59.4296\\
16	102.549\\
8	198.758\\
4	394.548\\
2	791.752\\
};
\end{axis}
\begin{axis}[%
width=0.36\columnwidth,
height=0.196\columnwidth,
at={(0\columnwidth,0\columnwidth)},
scale only axis,
xmode=log,
xmin=1,
xmax=512,
xtick={1,2,4,8,16,32,64,128,256,512},
xticklabels={},
xminorticks=true,
separate axis lines,
every outer y axis line/.append style={mycolor2},
every y tick label/.append style={font=\color{mycolor2}},
every y tick/.append style={mycolor2},
ymin=13,
ymax=16,
yminorticks=true,
ytick={13,14,15,16},
axis background/.style={fill=none},
yticklabel pos=right,
]
\addplot [color=mycolor2, dashed, line width=2.0pt, forget plot]
  table[row sep=crcr]{%
512	15\\
256	14\\
128	15\\
64	14\\
32	14\\
16	14\\
8	14\\
4	14\\
2	14\\
};
\end{axis}

\begin{axis}[%
width=0.36\columnwidth,
height=0.196\columnwidth,
at={(0.44\columnwidth,0\columnwidth)},
scale only axis,
xmode=log,
xmin=1,
xmax=512,
xtick={  1,   2,   4,   8,  16,  32,  64, 128, 256, 512},
xticklabels={  1,   2,   4,   8,  16,  32,  64, 128, 256, 512},
x tick label style={rotate=45,font=\small},
xminorticks=true,
xlabel style={font=\color{mycolor3}},
xlabel={MPI Tasks},
separate axis lines,
every outer y axis line/.append style={mycolor1},
every y tick label/.append style={font=\color{mycolor1}\small},
every y tick/.append style={mycolor1},
ymode=log,
ymin=3,
ymax=470,
ymode=log,
axis background/.style={fill=none},
title style={font=\small\color{mycolor3}},
title={$p = 5$},
title style={yshift=-1em},
yticklabel pos=left,
xmajorgrids,
xminorgrids,
ymajorgrids,
yminorgrids
]
\addplot [color=mycolor1, line width=2.0pt, mark=o, mark options={solid, mycolor1}, forget plot]
  table[row sep=crcr]{%
512	3.32858\\
256	3.47113\\
128	5.03969\\
64	9.55767\\
32	16.6136\\
16	30.1514\\
8	58.8007\\
4	113.04\\
2	227.898\\
1	436.502\\
};
\end{axis}
\begin{axis}[%
width=0.36\columnwidth,
height=0.196\columnwidth,
at={(0.44\columnwidth,0\columnwidth)},
scale only axis,
xmode=log,
xmin=1,
xmax=512,
xtick={  1,   2,   4,   8,  16,  32,  64, 128, 256, 512},
xticklabels={  1,   2,   4,   8,  16,  32,  64, 128, 256, 512},
x tick label style={rotate=45,font=\small},
xminorticks=true,
separate axis lines,
every outer y axis line/.append style={mycolor2},
every y tick label/.append style={font=\color{mycolor2}\small},
every y tick/.append style={mycolor2},
ymin=14,
ymax=16,
axis background/.style={fill=none},
ytick={14,15,16},
ylabel style={font=\color{mycolor2}},
ylabel={Iterations},
yticklabel pos=right,
xmajorgrids,
xminorgrids,
ymajorgrids,
yminorgrids
]
\addplot [color=mycolor2, dashed, line width=2.0pt, forget plot]
  table[row sep=crcr]{%
512	15\\
256	15\\
128	14\\
64	15\\
32	14\\
16	15\\
8	15\\
4	15\\
2	15\\
1	15\\
};
\end{axis}

\begin{axis}[%
width=0.36\columnwidth,
height=0.196\columnwidth,
at={(0\columnwidth,-0.27\columnwidth)},
scale only axis,
xmode=log,
xmin=1,
xmax=512,
xtick={1,2,4,8,16,32,64,128,256,512},
xticklabels={1,2,4,8,16,32,64,128,256,512},
xticklabel style={rotate=45,font=\small,anchor=east},
xminorticks=true,
xlabel style={font=\color{mycolor3}},
xlabel={MPI Tasks},
separate axis lines,
every outer y axis line/.append style={mycolor1},
every y tick label/.append style={font=\color{mycolor1}},
every y tick/.append style={mycolor1},
ymode=log,
ymin=5,
ymax=900,
yminorticks=true,
ylabel style={font=\color{mycolor1}},
ylabel={Time (s)},
axis background/.style={fill=white},
title style={font=\small},
title={$p = 6$},
title style={yshift=-1em},
yticklabel pos=left,
xmajorgrids,
xminorgrids,
ymajorgrids,
yminorgrids
]
\addplot [color=mycolor1, line width=2.0pt, mark=o, mark options={solid, mycolor1}, forget plot]
  table[row sep=crcr]{%
512	5.27032\\
256	5.99651\\
128	9.9437\\
64	18.1747\\
32	33.4493\\
16	56.9403\\
8	119.058\\
4	228.692\\
2	454.154\\
1	874.468\\
};
\end{axis}
\begin{axis}[%
width=0.36\columnwidth,
height=0.196\columnwidth,
at={(0\columnwidth,-0.27\columnwidth)},
scale only axis,
xmode=log,
xmin=1,
xmax=512,
xtick={1,2,4,8,16,32,64,128,256,512},
xticklabels={},
xminorticks=true,
separate axis lines,
every outer y axis line/.append style={mycolor2},
every y tick label/.append style={font=\color{mycolor2}},
every y tick/.append style={mycolor2},
ymin=14,
ymax=17,
yminorticks=true,
ytick={15,16,17},
axis background/.style={fill=none},
yticklabel pos=right,
ylabel style={font=\color{mycolor2}},
ylabel={Iterations},
]
\addplot [color=mycolor2, dashed, line width=2.0pt, forget plot]
  table[row sep=crcr]{%
512	16\\
256	15\\
128	15\\
64	15\\
32	15\\
16	15\\
8	16\\
4	16\\
2	16\\
1	16\\
};
\end{axis}

\end{tikzpicture}%
    \caption{Strong scaling for the 3D $L$-shaped domain problem for a number of task running from $1$ to $512$ {on the largest available test case for each $p$, i.e., sizes $2691168$, $2775752$, $385050$, $408408$, and $432692$ respectively.}}
    \label{fig:lshaped_h_strong_scaling}
\end{figure}

Figure~\ref{fig:lshaped_h_strong_scaling} presents the corresponding time-to-solution (left axis) and iteration counts (right axis). We observe substantial reductions in wall-clock time across the entire scaling range and nearly constant iteration counts in most configurations, indicating that the preconditioner’s effectiveness is largely preserved under distribution. Only at the highest task counts and for larger values of $p$ do we observe a mild increase in the number of iterations; however, the overall time-to-solution continues to decrease, showing that performance remains dominated by computational workload rather than convergence degradation.
Quantifying the parallel performance, for $p=2$ the wall-clock time decreases from $87.09$,s on a single task to $0.69$,s on $512$ tasks, giving a speedup of $125.5\times$ and a parallel efficiency of $24.5\%$. For $p=3$, the time drops from $359.35$,s to $1.94$,s, corresponding to a speedup of $184.9\times$ and an efficiency of $36.1\%$. For $p=4$, a single-task run could not be completed due to memory constraints; taking the two-task time ($791.75$,s) as baseline, the reduction to $6.76$,s at $512$ tasks yields a relative speedup of $117.2\times$ with an efficiency of $45.8\%$ (relative to the 256-fold increase in task count). The $p=5$ case scales from $436.50$,s to $3.33$,s (speedup $131.2\times$, efficiency $25.6\%$), while for $p=6$ the time decreases from $874.47$,s to $5.27$,s (speedup $166.0\times$, efficiency $32.4\%$). Across all degrees, iteration counts remain essentially constant—approximately $10$ for $p=2$ and $p=3$, $14$--$15$ for $p=4$ and $p=5$, and $15$--$16$ for $p=6$—indicating that the preconditioner maintains its algorithmic robustness on the multipatch geometry. Although efficiencies at 512 tasks are moderate for the lower degrees, the absolute speedups confirm that the presence of patch interfaces does not impede scalability; the higher efficiency observed for $p=4$ relative to its baseline further reflects favorable workload distribution once a sufficiently large local problem size per task is reached.

\subsection{Non-isoparametric test case}\label{sec:noniso}
Finally, we consider the non-isoparametric configuration introduced in Section~\ref{sec:iso-or-noniso}: the geometry is represented with NURBS while the solution is discretized with B-splines on the quarter-of-ring domain. Within each patch, we use maximum regularity ($\mathcal{C}^{p-1}$). The solver settings are unchanged with respect to the previous studies. %
Table~\ref{tab:poisson_ring} reports single-process results as $k$ and $p$ grow. Iteration counts remain modest and only mildly dependent on $k$, while operator complexity increases with $p$ in a controlled manner. Overall, the non-isoparametric discretization exhibits convergence behavior comparable to the cube and multipatch $L$-shaped cases, indicating that the change of geometric mapping basis does not adversely affect the preconditioner’s effectiveness.
\begin{table}[htbp]
    \centering
    \begin{subtable}{0.2\columnwidth}
    \begin{tabular}{l|cccc}
    \toprule
    $k \backslash p$  &  2 & 3 & 4 & 5\\
    \midrule
12              & 8               & 10              & 12              & 20             \\
24              & 8               & 10              & 16              & 24             \\
48              & 8               & 11              & 16              & 28             \\
96              & 9               & 11              & 16              & 28             \\
    \bottomrule
    \end{tabular}
    \caption{Iteration count}
      \end{subtable}\hfil
    \begin{subtable}{0.42\columnwidth}
    \begin{tabular}{l|cccc}
    \toprule
    $k \backslash p$  &  2 & 3 & 4 & 5\\
    \midrule
12              & 2184            & 2730            & 3360            & 4080            \\
24              & 15600           & 17550           & 19656           & 21924           \\
48              & 117600          & 124950          & 132600          & 140556          \\
96              & 912576          & 941094          & 970200          & 999900          \\   \bottomrule
    \end{tabular}
    \caption{Matrix size}
    \end{subtable}\hfil
    \begin{subtable}{0.25\columnwidth}
    \begin{tabular}{l|cccc}
    \toprule
    $k \backslash p$  &  2 & 3 & 4 & 5\\
    \midrule

12              & 1.22 & 1.15 & 1.10 & 1.08 \\
24              & 1.34 & 1.23 & 1.18 & 1.14 \\
48              & 1.44 & 1.29 & 1.25 & 1.18 \\
96              & 1.49 & 1.34 & 1.29 & 1.21 \\
    \bottomrule
    \end{tabular}
    \caption{Operator complexity}
    \end{subtable}
  
        \caption{Result for the Poisson on the quarter of ring domain case on a single task for varying number of subdivision $k$ of the base mesh, and degree $p$ of the B-Spline basis.}
    \label{tab:poisson_ring}
\end{table}
\begin{figure}[htbp]
    \centering
    \includegraphics[width=\columnwidth]{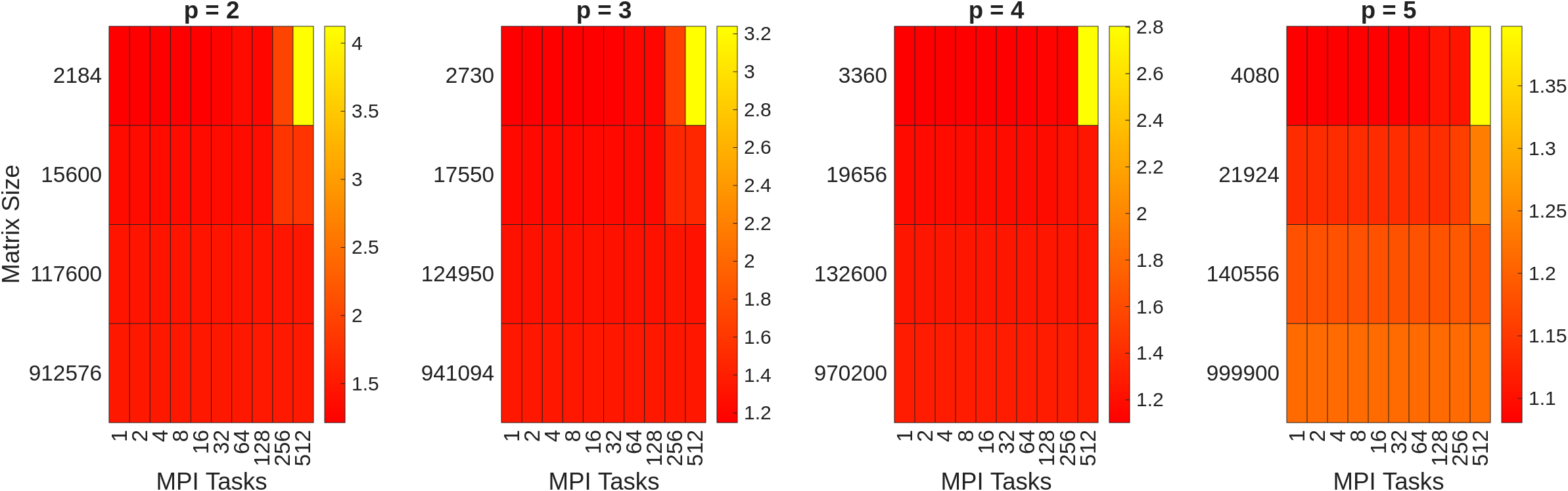}
    \caption{Operator complexity for the non-isoparametric quarter of ring case for increasing values of $k$ and $p$, and number of MPI tasks running from $1$ to $512$.}
    \label{fig:ring_opc}
\end{figure}
As in the other test cases, distributing the graph across many MPI tasks slightly increases the operator complexity for the smallest matrix sizes (Fig.~\ref{fig:ring_opc}), reflecting a less effective matching at very fine partitions. For medium and large problems, the complexity remains close to the single-task values across the entire range of task counts.
\begin{figure}[hbt]
    \definecolor{mycolor1}{rgb}{0.06667,0.44314,0.74510}%
\definecolor{mycolor2}{rgb}{0.86667,0.32941,0.00000}%
\definecolor{mycolor3}{rgb}{0.12941,0.12941,0.12941}%
\begin{tikzpicture}

\begin{axis}[%
width=0.36\columnwidth,
height=0.196\columnwidth,
at={(0\columnwidth,0.27\columnwidth)},
scale only axis,
xmode=log,
xmin=1,
xmax=512,
xtick={1,2,4,8,16,32,64,128,256,512},
xticklabels={1,2,4,8,16,32,64,128,256,512},
xticklabel style={rotate=45,font=\footnotesize,anchor=east},
xminorticks=true,
xlabel style={font=\color{mycolor3}},
separate axis lines,
every outer y axis line/.append style={mycolor1},
every y tick label/.append style={font=\color{mycolor1}},
every y tick/.append style={mycolor1},
ymode=log,
ymin=0.35,
ymax=30,
yminorticks=true,
ylabel style={font=\color{mycolor1}},
ylabel={Time (s)},
axis background/.style={fill=white},
title style={font=\small},
title={$p = 2$},
title style={yshift=-1em},
yticklabel pos=left,
xmajorgrids,
xminorgrids,
ymajorgrids,
yminorgrids
]
\addplot [color=mycolor1, line width=2.0pt, mark=o, mark options={solid, mycolor1}, forget plot]
  table[row sep=crcr]{%
512	0.489061\\
256	0.379437\\
128	0.37472\\
64	0.553255\\
32	1.0342\\
16	1.7188\\
8	3.3278\\
4	6.27868\\
2	12.4676\\
1	24.5423\\
};
\end{axis}
\begin{axis}[%
width=0.36\columnwidth,
height=0.196\columnwidth,
at={(0\columnwidth,0.27\columnwidth)},
scale only axis,
xmode=log,
xmin=1,
xmax=512,
xtick={1,2,4,8,16,32,64,128,256,512},
xticklabels={},
xminorticks=true,
xlabel style={font=\color{mycolor3}},
separate axis lines,
every outer y axis line/.append style={mycolor2},
every y tick label/.append style={font=\color{mycolor2}},
every y tick/.append style={mycolor2},
ymin=8.5,
ymax=9.5,
ytick={9},
yticklabels={9},
yminorticks=true,
axis background/.style={fill=none},
title style={font=\bfseries\color{mycolor3}},
yticklabel pos=right,
]
\addplot [color=mycolor2, dashed, line width=2.0pt, forget plot]
  table[row sep=crcr]{%
512	9\\
256	9\\
128	9\\
64	9\\
32	9\\
16	9\\
8	9\\
4	9\\
2	9\\
1	9\\
};
\end{axis}

\begin{axis}[%
width=0.36\columnwidth,
height=0.196\columnwidth,
at={(0.44\columnwidth,0.27\columnwidth)},
scale only axis,
xmode=log,
xmin=1,
xmax=512,
xtick={1,2,4,8,16,32,64,128,256,512},
xticklabels={1,2,4,8,16,32,64,128,256,512},
xticklabel style={rotate=45,font=\small,anchor=east},
xminorticks=true,
xlabel style={font=\color{mycolor3}},
separate axis lines,
every outer y axis line/.append style={mycolor1},
every y tick label/.append style={font=\color{mycolor1}},
every y tick/.append style={mycolor1},
ymode=log,
ymin=1.0,
ymax=140,
ymode=log,
yminorticks=true,
ylabel style={font=\color{mycolor1}},
axis background/.style={fill=white},
title style={font=\small},
title={$p = 3$},
title style={yshift=-1em},
yticklabel pos=left,
xmajorgrids,
xminorgrids,
ymajorgrids,
yminorgrids
]
\addplot [color=mycolor1, line width=2.0pt, mark=o, mark options={solid, mycolor1}, forget plot]
  table[row sep=crcr]{%
512	1.16442\\
256	1.19635\\
128	1.55532\\
64	2.80519\\
32	5.35713\\
16	9.16182\\
8	17.6221\\
4	33.9123\\
2	66.6179\\
1	131.453\\
};
\end{axis}
\begin{axis}[%
width=0.36\columnwidth,
height=0.196\columnwidth,
at={(0.44\columnwidth,0.27\columnwidth)},
scale only axis,
xmode=log,
xmin=1,
xmax=512,
xtick={1,2,4,8,16,32,64,128,256,512},
xticklabels={},
xminorticks=true,
separate axis lines,
every outer y axis line/.append style={mycolor2},
every y tick label/.append style={font=\color{mycolor2}},
every y tick/.append style={mycolor2},
ymin=10.5,
ymax=11.5,
ytick={11},
yticklabels={11},
yminorticks=true,
axis background/.style={fill=none},
yticklabel pos=right,
ylabel style={font=\color{mycolor2}},
ylabel={Iterations},
]
\addplot [color=mycolor2, dashed, line width=2.0pt, forget plot]
  table[row sep=crcr]{%
512	11\\
256	11\\
128	11\\
64	11\\
32	11\\
16	11\\
8	11\\
4	11\\
2	11\\
1	11\\
};
\end{axis}

\begin{axis}[%
width=0.36\columnwidth,
height=0.196\columnwidth,
at={(0\columnwidth,0\columnwidth)},
scale only axis,
xmode=log,
xmin=1,
xmax=512,
xtick={1,2,4,8,16,32,64,128,256,512},
xticklabels={1,2,4,8,16,32,64,128,256,512},
xticklabel style={rotate=45,font=\small,anchor=east},
xminorticks=true,
xlabel style={font=\color{mycolor3}},
xlabel={MPI Tasks},
separate axis lines,
every outer y axis line/.append style={mycolor1},
every y tick label/.append style={font=\color{mycolor1}},
every y tick/.append style={mycolor1},
ymode=log,
ymin=3,
ymax=600,
yminorticks=true,
ylabel style={font=\color{mycolor1}},
ylabel={Time (s)},
axis background/.style={fill=white},
title style={font=\small},
title={$p = 4$},
title style={yshift=-1em},
yticklabel pos=left,
xmajorgrids,
xminorgrids,
ymajorgrids,
yminorgrids
]
\addplot [color=mycolor1, line width=2.0pt, mark=o, mark options={solid, mycolor1}, forget plot]
  table[row sep=crcr]{%
512	3.39318\\
256	4.57579\\
128	6.58098\\
64	12.2969\\
32	25.061\\
16	40.3607\\
8	77.8622\\
4	151.282\\
2	295.48\\
1	582.9\\
};
\end{axis}
\begin{axis}[%
width=0.36\columnwidth,
height=0.196\columnwidth,
at={(0\columnwidth,0\columnwidth)},
scale only axis,
xmode=log,
xmin=1,
xmax=512,
xtick={1,2,4,8,16,32,64,128,256,512},
xticklabels={},
xminorticks=true,
separate axis lines,
every outer y axis line/.append style={mycolor2},
every y tick label/.append style={font=\color{mycolor2}},
every y tick/.append style={mycolor2},
ymin=15.5,
ymax=17.5,
ytick={16,17},
yticklabels={16,17},
yminorticks=true,
axis background/.style={fill=none},
yticklabel pos=right,
]
\addplot [color=mycolor2, dashed, line width=2.0pt, forget plot]
  table[row sep=crcr]{%
512	16\\
256	17\\
128	16\\
64	16\\
32	17\\
16	16\\
8	16\\
4	16\\
2	16\\
1	16\\
};
\end{axis}

\begin{axis}[%
width=0.36\columnwidth,
height=0.196\columnwidth,
at={(0.44\columnwidth,0\columnwidth)},
scale only axis,
xmode=log,
xmin=1,
xmax=512,
xtick={  1,   2,   4,   8,  16,  32,  64, 128, 256, 512},
xticklabels={  1,   2,   4,   8,  16,  32,  64, 128, 256, 512},
x tick label style={rotate=45,font=\small},
xminorticks=true,
xlabel style={font=\color{mycolor3}},
xlabel={MPI Tasks},
separate axis lines,
every outer y axis line/.append style={mycolor1},
every y tick label/.append style={font=\color{mycolor1}\small},
every y tick/.append style={mycolor1},
ymode=log,
ymin=9,
ymax=2100,
ymode=log,
axis background/.style={fill=none},
title style={font=\small\color{mycolor3}},
title={$p = 5$},
title style={yshift=-1em},
yticklabel pos=left,
xmajorgrids,
xminorgrids,
ymajorgrids,
yminorgrids
]
\addplot [color=mycolor1, line width=2.0pt, mark=o, mark options={solid, mycolor1}, forget plot]
  table[row sep=crcr]{%
512	9.47577\\
256	13.0873\\
128	22.6291\\
64	42.2747\\
32	82.9075\\
16	137.004\\
8	272.378\\
4	529.73\\
2	1047.9\\
1	2032.64\\
};
\end{axis}
\begin{axis}[%
width=0.36\columnwidth,
height=0.196\columnwidth,
at={(0.44\columnwidth,0\columnwidth)},
scale only axis,
xmode=log,
xmin=1,
xmax=512,
xtick={  1,   2,   4,   8,  16,  32,  64, 128, 256, 512},
xticklabels={  1,   2,   4,   8,  16,  32,  64, 128, 256, 512},
x tick label style={rotate=45,font=\small},
xminorticks=true,
xlabel style={font=\color{mycolor3}},
xlabel={MPI Tasks},
separate axis lines,
every outer y axis line/.append style={mycolor2},
every y tick label/.append style={font=\color{mycolor2}\small},
every y tick/.append style={mycolor2},
ymin=26.5,
ymax=28.5,
ytick={27,28},
yticklabels={27,28},
axis background/.style={fill=none},
ylabel style={font=\color{mycolor2}},
ylabel={Iterations},
yticklabel pos=right,
xmajorgrids,
xminorgrids,
ymajorgrids,
yminorgrids
]
\addplot [color=mycolor2, dashed, line width=2.0pt, forget plot]
  table[row sep=crcr]{%
512	28\\
256	28\\
128	28\\
64	28\\
32	28\\
16	27\\
8	28\\
4	28\\
2	28\\
1	28\\
};
\end{axis}

\end{tikzpicture}%
    \caption{Strong scaling for the non-isoparametric quarter of ring problem for a number of task running from $1$ to $512$ {on the largest available test case for each $p$, i.e., sizes $912576$, $941094$, $970200$ and $999900$ respectively.}}
    \label{fig:ring_h_strong_scaling}
\end{figure}
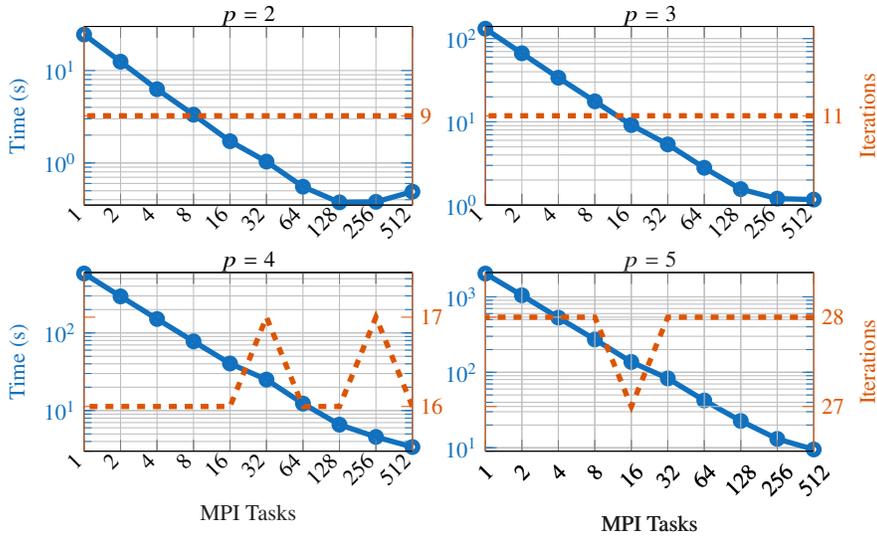
The strong-scaling results are reported in Fig.~\ref{fig:ring_h_strong_scaling}. For $p=2$, the wall-clock time decreases from $24.54\,\mathrm{s}$ on one task to $0.49\,\mathrm{s}$ on $512$ tasks, giving a speedup of $50.2\times$ (parallel efficiency $9.8\%$); the modest efficiency at the largest task count reflects the small local problem size per rank. For $p=3$, the time drops from $131.45\,\mathrm{s}$ to $1.16\,\mathrm{s}$, corresponding to a speedup of $112.9\times$ and an efficiency of $22.0\%$. The $p=4$ case improves from $582.90\,\mathrm{s}$ to $3.39\,\mathrm{s}$ (speedup $171.8\times$, efficiency $33.6\%$), while for $p=5$ the time decreases from $2032.64\,\mathrm{s}$ to $9.48\,\mathrm{s}$ (speedup $214.5\times$, efficiency $41.9\%$). These results show a steady increase in parallel efficiency with higher polynomial degree, attributable to the larger per-task workload mitigating communication and setup overheads. Iteration counts remain nearly constant across the scaling range---about $9$ for $p=2$, $11$ for $p=3$, $16$--$17$ for $p=4$, and $27$--$28$ for $p=5$---indicating that convergence robustness is preserved despite the non-isoparametric geometry and the extreme task counts. Overall, the non-isoparametric configuration attains strong-scaling performance comparable to, and for higher $p$ even slightly better than, the multipatch and cube cases, confirming that differences in geometry and solution spaces do not introduce detrimental overhead in the parallel solver.

\subsection{Experiments with GPU support}
In this section, we focus on the performance attainable when using GPUs to execute the solve phase of the different linear systems. As discussed in Section~\ref{sec:psctoolkit}, from the user’s perspective this requires nothing more than setting the appropriate \lstinline[language=Fortran,style=fortran]|mold=| options for matrices, vectors, and descriptors. Nowadays, GPUs suitable for general-purpose computation are also available on laptops, workstations, and small servers. Since \texttt{PSCToolkit} can be compiled and executed in such environments as well, we include performance results for these systems in Section~\ref{sec:small_systems} to provide a more complete overview.

In this part of the analysis, we employ the hacked-\emph{Ellpack} format~\cite{10.1145/3017994} as the matrix storage scheme on the GPU, which is represented in the library by the string \lstinline[language=Fortran,style=fortran]|'HLG'|. In a cluster environment, the number of available GPUs per node is typically smaller than the number of CPU cores that could be assigned to MPI tasks. Although MPI ranks do not strictly need to match the number of GPUs, running more than one rank per GPU generally leads to contention and reduced performance. Therefore, since each node in our configuration is equipped with four GPUs, we restrict execution to at most four MPI tasks per node, adopting a one-rank-per-GPU strategy.

\begin{figure}[htb]
    \centering
    \begin{subfigure}[t]{0.3\columnwidth}
    \definecolor{mycolor1}{rgb}{0.06600,0.44300,0.74500}%
\definecolor{mycolor2}{RGB}{118,185,0}%
\definecolor{mycolor3}{rgb}{0.12941,0.12941,0.12941}%
\begin{tikzpicture}

\begin{axis}[%
width=0.8\columnwidth,
height=0.8\columnwidth,
at={(0\columnwidth,0\columnwidth)},
scale only axis,
xmode=log,
xmin=1,
xmax=512,
xtick={  1,   2,   4,   8,  16,  32,  64, 128, 256, 512},
xticklabels={  1,   2,   4,   8,  16,  32,  64, 128, 256, 512},
xticklabel style={rotate=45,font=\footnotesize,anchor=east},
xminorticks=true,
xlabel style={font=\color{mycolor3}},
xlabel={MPI Tasks},
ymode=log,
ymin=0.261251,
ymax=65.7121,
yminorticks=true,
ylabel style={font=\color{mycolor3}},
ylabel={Time (s)},
axis background/.style={fill=white},
xmajorgrids,
xminorgrids,
ymajorgrids,
yminorgrids,
legend columns=2,
legend style={legend cell align=left, align=left, draw=none, fill=none, font=\scriptsize, at={(0.5,1.05)}, anchor=south}
]
\addplot [color=mycolor1, line width=2.0pt, mark=o, mark options={solid, mycolor1}]
  table[row sep=crcr]{%
512	0.261251\\
256	0.543746\\
128	1.06254\\
64	4.81789\\
32	5.57041\\
16	6.80025\\
8	10.0268\\
4	17.1933\\
2	37.072\\
1	65.7121\\
};
\addlegendentry{MPI}

\addplot [color=mycolor2, dashed, line width=2.0pt, mark=square, mark options={solid, mycolor2}]
  table[row sep=crcr]{%
64	0.392815\\
32	0.360067\\
16	0.370575\\
8	0.455907\\
4	0.557315\\
2	0.88264\\
1	1.56587\\
};
\addlegendentry{GPU - A30}

\node[above left, align=right, inner sep=0, yshift=1em, font=\color{mycolor3}]
at (axis cs:512,0.261) {16};
\node[above left, align=right, inner sep=0, yshift=1em, font=\color{mycolor3}]
at (axis cs:256,0.544) {8};
\node[above left, align=right, inner sep=0,yshift=1em, font=\color{mycolor3}]
at (axis cs:128,1.063) {4};
\node[above left, align=right, inner sep=0, yshift=1em, font=\color{mycolor3}]
at (axis cs:64,4.818) {2};
\node[above left, align=right, inner sep=0, yshift=1em, font=\color{mycolor3}]
at (axis cs:32,5.57) {1};
\node[above left, align=right, inner sep=0, yshift=1em, font=\color{mycolor3}]
at (axis cs:16,6.8) {1};
\node[above left, align=right, inner sep=0, yshift=1em, font=\color{mycolor3}]
at (axis cs:8,10.027) {1};
\node[above left, align=right, inner sep=0, yshift=-1em, font=\color{mycolor3}]
at (axis cs:4,17.193) {1};
\node[above left, align=right, inner sep=0, yshift=-1em, font=\color{mycolor3}]
at (axis cs:2,37.072) {1};
\node[above left, align=right, inner sep=0, yshift=-1em, font=\color{mycolor3}]
at (axis cs:1,65.712) {1};
\node[above left, align=right, inner sep=0, yshift=1em, font=\color{mycolor3}]
at (axis cs:64,0.393) {16};
\node[above left, align=right, inner sep=0, yshift=1em, font=\color{mycolor3}]
at (axis cs:32,0.36) {8};
\node[above left, align=right, inner sep=0, yshift=1em, font=\color{mycolor3}]
at (axis cs:16,0.371) {4};
\node[above left, align=right, inner sep=0, yshift=1em, font=\color{mycolor3}]
at (axis cs:8,0.456) {2};
\node[above left, align=right, inner sep=0, yshift=1em, font=\color{mycolor3}]
at (axis cs:4,0.557) {1};
\node[above left, align=right, inner sep=0, yshift=1em, font=\color{mycolor3}]
at (axis cs:2,0.883) {1};
\node[above left, align=right, inner sep=0, yshift=1em, font=\color{mycolor3}]
at (axis cs:1,1.566) {1};
\end{axis}

\end{tikzpicture}%

    \caption{3D Cube, $k = 96$, $p = 3$}
    \end{subfigure}\hspace{1.5em}
    \begin{subfigure}[t]{0.3\columnwidth}
    \definecolor{mycolor1}{rgb}{0.06600,0.44300,0.74500}%
\definecolor{mycolor2}{RGB}{118,185,0}%
\definecolor{mycolor3}{rgb}{0.12941,0.12941,0.12941}%
\begin{tikzpicture}

\begin{axis}[%
width=0.8\columnwidth,
height=0.8\columnwidth,
at={(0\columnwidth,0\columnwidth)},
scale only axis,
xmode=log,
xmin=1,
xmax=512,
xtick={  1,   2,   4,   8,  16,  32,  64, 128, 256, 512},
xticklabels={  1,   2,   4,   8,  16,  32,  64, 128, 256, 512},
xticklabel style={rotate=45,font=\footnotesize,anchor=east},
xminorticks=true,
xlabel style={font=\color{mycolor3}},
xlabel={MPI Tasks},
ymode=log,
ymin=0.393781,
ymax=184.085,
yminorticks=true,
axis background/.style={fill=white},
xmajorgrids,
xminorgrids,
ymajorgrids,
yminorgrids,
legend columns=2,
legend style={legend cell align=left, align=left, draw=none, fill=none, font=\scriptsize, at={(0.5,1.05)}, anchor=south}
]
\addplot [color=mycolor1, line width=2.0pt, mark=o, mark options={solid, mycolor1}]
  table[row sep=crcr]{%
512	0.771403\\
256	1.98783\\
128	3.64505\\
64	5.41317\\
32	10.1141\\
16	15.6049\\
8	27.9861\\
4	48.0863\\
2	110.305\\
1	184.085\\
};
\addlegendentry{MPI}

\addplot [color=mycolor2, dashed, line width=2.0pt, mark=square, mark options={solid, mycolor2}]
  table[row sep=crcr]{%
64	0.393781\\
32	0.457148\\
16	0.59161\\
8	0.864615\\
4	1.36641\\
2	2.46736\\
1	4.63595\\
};
\addlegendentry{GPU - A30}

\node[above, yshift=3pt, align=center, inner sep=0, font=\color{mycolor3}]
at (axis cs:512,0.771403) {16};
\node[above, yshift=3pt, align=center, inner sep=0, font=\color{mycolor3}]
at (axis cs:256,1.98783) {8};
\node[above, yshift=3pt, align=center, inner sep=0, font=\color{mycolor3}]
at (axis cs:128,3.64505) {4};
\node[above, yshift=3pt, align=center, inner sep=0, font=\color{mycolor3}]
at (axis cs:64,5.41317) {2};
\node[above, yshift=3pt, align=center, inner sep=0, font=\color{mycolor3}]
at (axis cs:32,10.1141) {1};
\node[above, yshift=3pt, align=center, inner sep=0, font=\color{mycolor3}]
at (axis cs:16,15.6049) {1};
\node[above, yshift=3pt, align=center, inner sep=0, font=\color{mycolor3}]
at (axis cs:8,27.9861) {1};
\node[above, yshift=3pt, align=center, inner sep=0, font=\color{mycolor3}]
at (axis cs:4,48.0863) {1};
\node[above, yshift=3pt, align=center, inner sep=0, font=\color{mycolor3}]
at (axis cs:2,110.305) {1};
\node[above, yshift=3pt, align=center, inner sep=0, font=\color{mycolor3}]
at (axis cs:1,184.085) {1};
\node[above, yshift=3pt, align=center, inner sep=0, font=\color{mycolor3}]
at (axis cs:64,0.393781) {16};
\node[above, yshift=3pt, align=center, inner sep=0, font=\color{mycolor3}]
at (axis cs:32,0.457148) {8};
\node[above, yshift=3pt, align=center, inner sep=0, font=\color{mycolor3}]
at (axis cs:16,0.59161) {4};
\node[above, yshift=3pt, align=center, inner sep=0, font=\color{mycolor3}]
at (axis cs:8,0.864615) {2};
\node[above, yshift=3pt, align=center, inner sep=0, font=\color{mycolor3}]
at (axis cs:4,1.36641) {1};
\node[above, yshift=3pt, align=center, inner sep=0, font=\color{mycolor3}]
at (axis cs:2,2.46736) {1};
\node[above, yshift=3pt, align=center, inner sep=0, font=\color{mycolor3}]
at (axis cs:1,4.63595) {1};
\end{axis}
\end{tikzpicture}%

    \caption{Multipatch, $k=96$, $p=3$}
    \end{subfigure}
    \begin{subfigure}[t]{0.3\columnwidth}
    \definecolor{mycolor1}{rgb}{0.06600,0.44300,0.74500}%
\definecolor{mycolor2}{RGB}{118,185,0}%
\definecolor{mycolor3}{rgb}{0.12941,0.12941,0.12941}%
\begin{tikzpicture}

\begin{axis}[%
width=0.8\columnwidth,
height=0.8\columnwidth,
at={(0\columnwidth,0\columnwidth)},
scale only axis,
xmode=log,
xmin=1,
xmax=512,
xtick={  1,   2,   4,   8,  16,  32,  64, 128, 256, 512},
xticklabels={  1,   2,   4,   8,  16,  32,  64, 128, 256, 512},
xticklabel style={rotate=45,font=\footnotesize,anchor=east},
xminorticks=true,
xlabel style={font=\color{mycolor3}},
xlabel={MPI Tasks},
ymode=log,
ymin=0.971237,
ymax=307.065,
yminorticks=true,
axis background/.style={fill=white},
xmajorgrids,
xminorgrids,
ymajorgrids,
yminorgrids,
legend columns=2,
legend style={legend cell align=left, align=left, draw=none, fill=none, font=\scriptsize, at={(0.5,1.05)}, anchor=south}
]
\addplot [color=mycolor1, line width=2.0pt, mark=o, mark options={solid, mycolor1}]
  table[row sep=crcr]{%
512	1.14938\\
256	2.18353\\
128	4.48546\\
64	14.7348\\
32	18.1395\\
16	26.2885\\
8	66.2324\\
4	124.075\\
2	164.072\\
1	307.065\\
};
\addlegendentry{MPI}

\addplot [color=mycolor2, dashed, line width=2.0pt, mark=square, mark options={solid, mycolor2}]
  table[row sep=crcr]{%
64	0.971237\\
32	1.01694\\
16	1.07725\\
8	1.40635\\
4	2.01799\\
2	3.39436\\
1	6.23353\\
};
\addlegendentry{GPU - A30}

\node[above, align=right, inner sep=2pt, font=\color{mycolor3}]
at (axis cs:512,1.149) {16};
\node[above, align=right, inner sep=2pt, font=\color{mycolor3}]
at (axis cs:256,2.184) {8};
\node[above, align=right, inner sep=2pt, font=\color{mycolor3}]
at (axis cs:128,4.485) {4};
\node[above, align=right, inner sep=2pt, font=\color{mycolor3}]
at (axis cs:64,14.735) {2};
\node[above, align=right, inner sep=2pt, font=\color{mycolor3}]
at (axis cs:32,18.139) {1};
\node[above, align=right, inner sep=2pt, font=\color{mycolor3}]
at (axis cs:16,26.288) {1};
\node[above, align=right, inner sep=2pt, font=\color{mycolor3}]
at (axis cs:8,66.232) {1};
\node[above, align=right, inner sep=2pt, font=\color{mycolor3}]
at (axis cs:4,124.075) {1};
\node[above, align=right, inner sep=2pt, font=\color{mycolor3}]
at (axis cs:2,164.072) {1};
\node[above, align=right, inner sep=2pt, font=\color{mycolor3}]
at (axis cs:1,307.065) {1};
\node[above, align=right, inner sep=2pt, font=\color{mycolor3}]
at (axis cs:64,0.971) {16};
\node[above, align=right, inner sep=2pt, font=\color{mycolor3}]
at (axis cs:32,1.017) {8};
\node[above, align=right, inner sep=2pt, font=\color{mycolor3}]
at (axis cs:16,1.077) {4};
\node[above, align=right, inner sep=2pt, font=\color{mycolor3}]
at (axis cs:8,1.406) {2};
\node[above, align=right, inner sep=2pt, font=\color{mycolor3}]
at (axis cs:4,2.018) {1};
\node[above, align=right, inner sep=2pt, font=\color{mycolor3}]
at (axis cs:2,3.394) {1};
\node[above, align=right, inner sep=2pt, font=\color{mycolor3}]
at (axis cs:1,6.234) {1};
\end{axis}
\end{tikzpicture}%

    \caption{Ring, $k = 96$, $p = 4$}
    \end{subfigure}%
    
    \caption{Strong scaling for the GPU. The nodes of the Amelia cluster are equipped with 4 NVIDIA A30 GPUs per node, hence we can employ at most 4 MPI tasks per node. This means that using up to 16 nodes corresponds to at most 512 MPI tasks (32 tasks per node), and 64 GPUs (4 tasks per node).}
    \label{fig:gpu_strong_scaling}
\end{figure}
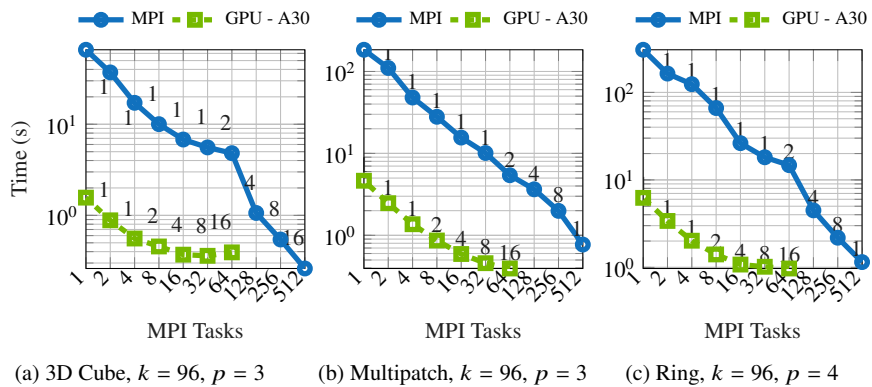
We report in Fig.~\ref{fig:gpu_strong_scaling} a comparison between the strong-scaling solve time employing a number of GPUs running from 1 to 64, i.e., a number of computing nodes running from 1 to 16, and the solve time using the pure MPI implementation discussed in detail in Sections~\ref{sec:cube}, \ref{sec:multipatch-experiments} and~\ref{sec:noniso}.
The GPU implementation demonstrates substantial performance improvements across all three test cases. For the 3D cube problem with $k=96$ and $p=3$, a single GPU achieves a solve time of $1.57\,\mathrm{s}$ compared to $65.71\,\mathrm{s}$ for a single MPI task, yielding a speedup of approximately $42\times$. When scaling to 16 GPUs, the time reduces to $0.39\,\mathrm{s}$, outperforming even the 64-task MPI configuration ($4.82\,\mathrm{s}$) by a factor of about $12.3\times$. For the multipatch L-shaped domain with $k=96$ and $p=3$, the acceleration is similarly impressive: a single GPU solves the problem in $4.64\,\mathrm{s}$ versus $184.09\,\mathrm{s}$ for one MPI task (speedup $\approx 39.7\times$), while 16 GPUs require only $0.39\,\mathrm{s}$ compared to $5.41\,\mathrm{s}$ for 64 MPI tasks (speedup $\approx 13.7\times$). The non-isoparametric ring configuration with $k=96$ and $p=4$ exhibits the most dramatic gains: a single GPU completes the solve in $6.23\,\mathrm{s}$ versus $307.07\,\mathrm{s}$ on one MPI task (speedup $\approx 49.3\times$), and 16 GPUs achieve $0.97\,\mathrm{s}$ compared to $14.73\,\mathrm{s}$ for 64 MPI tasks (speedup $\approx 15.2\times$).

However, the performance advantage of GPUs diminishes when comparing 64 GPUs against 512 MPI tasks. At this extreme level of parallelization, each GPU processes a very small local matrix block. For the cube test case, 64 GPUs achieve $0.39\,\mathrm{s}$ compared to $0.26\,\mathrm{s}$ for 512 MPI tasks, indicating that the pure MPI configuration is actually faster. For the multipatch and ring cases, 64 GPUs achieve solve times of $0.39\,\mathrm{s}$ and $0.97\,\mathrm{s}$ versus $0.77\,\mathrm{s}$ and $1.15\,\mathrm{s}$ for 512 MPI tasks, corresponding to modest speedups of $2\times$ and $1.2\times$. This reduction in GPU efficiency at very high task counts is due primarily to two factors. First, as the problem is partitioned across more processing units, the local matrices become smaller, reducing computational intensity and limiting the GPUs' ability to exploit their high throughput and massive parallelism. Second, when local workloads become too small, data-movement overheads and kernel-launch latencies become comparable to the actual computation time, eroding the GPU advantage. In such fine-grained decompositions, the lightweight per-core overhead of MPI tasks may therefore become competitive with, or even superior to, GPU execution.

These results confirm that GPU acceleration provides significant computational advantages, particularly for problems with larger matrix sizes and higher polynomial degrees, where the dense computational kernels inherent to IgA discretizations are well suited to GPU architectures.

\subsubsection{Acceleration on Consumer GPUs}\label{sec:small_systems}
While our results thus far have emphasized the performance gains achieved on distributed systems within HPC environments, the \texttt{PSCToolkit} library is designed with portability and versatility in mind. In particular, it also enables substantial acceleration on consumer-grade hardware by exploiting GPUs commonly available in modern workstations, laptops, and mid-range servers. This capability broadens the accessibility of advanced numerical solvers, allowing users outside traditional HPC settings to benefit from significant performance improvements.
\begin{table}[htbp]
    \centering
    \begin{tabular}{lp{3.5cm}p{5.7cm}}
       \toprule
       \textbf{Machine}  & \textbf{CPU and GPU} & \textbf{Software environment}  \\
         \midrule
       Laptop  & Intel\textsuperscript{\textregistered} Core\textsuperscript{\texttrademark} i9-14900HX,\phantom{aa} NVIDIA GeForce RTX 4060 & \texttt{gcc} v13.3.0, \texttt{cuda} v12.6.2, \texttt{openmpi} v4.1.7, \texttt{openblas} v0.3.28  \\ 
       \midrule
       Workstation~1 & AMD Ryzen 5 PRO 5650G, NVIDIA T1000 8\,GB & \texttt{gcc} v12.2.0,  \texttt{cuda} v12.3.0, \texttt{openmpi} v5.0.2, \texttt{openblas} v0.3.26 \\
       \midrule
       Workstation~2 & Intel\textsuperscript{\textregistered} Core\textsuperscript{\texttrademark} i7-12700,\phantom{aaaaa} NVIDIA T1000 4\,GB & \texttt{llvm} v21.1.0, \texttt{cuda} v13.0, \texttt{mpich} v4.3.0, \texttt{openblas} v0.3.28 \\
       \midrule
       Server & AMD EPYC 7763,\phantom{NVIDIA} NVIDIA A40 & \texttt{gcc} v14.2.0, \texttt{cuda} v12.8.0, \texttt{openmpi} v4.1.8, \texttt{openblas} v0.3.28\\
         \bottomrule
    \end{tabular}
    \caption{Hardware and software configurations used to assess single-GPU acceleration.}
    \label{tab:consumer_hardware}
\end{table}

The experiments summarized in Table~\ref{tab:consumer_hardware} were carried out on systems featuring diverse combinations of hardware and software environments. This diversity reflects realistic user configurations, including variations in compiler toolchains (\texttt{gcc} and \texttt{llvm}), MPI implementations (\texttt{openmpi} and \texttt{mpich}), and CUDA toolkit versions ranging from 12.3 to 13.0. Despite these differences, the results demonstrate the robustness and portability of \texttt{PSCToolkit}, which can seamlessly adapt to heterogeneous environments without requiring any platform-specific adjustments or tuning. Such flexibility ensures that the same computational kernels and solver configurations can be deployed across different architectures and software stacks, further enhancing the library’s accessibility and reproducibility.

We evaluate the acceleration achieved on the systems and environments listed in Table~\ref{tab:consumer_hardware}, using the same benchmark problems discussed in Sections~\ref{sec:cube}, \ref{sec:multipatch-experiments}, and~\ref{sec:noniso}. Having already analyzed the algorithmic performance in terms of convergence and scalability, we now focus on the solve times and the resulting speedups obtained when comparing GPU execution against the corresponding CPU runs. 

\begin{figure}[htbp]
    \centering

    \begin{subfigure}{0.9\columnwidth}
\definecolor{nvidiaGreen}{RGB}{118,185,0}
\begin{tikzpicture}
\begin{axis}[
    ybar,
    bar width=12pt,
    width=0.9\columnwidth,
    height=4cm,
    enlarge x limits={0.2},
    ylabel={Time to solve (s)},
    symbolic x coords={Laptop,Workstation 1,Server,Workstation 2},
    xtick=data,
    legend style={
        at={(1.02,0.5)},
        anchor=west
    },
    nodes near coords,
    nodes near coords align={vertical},
    ymajorgrids=true
]

\addplot+[fill=blue!50] coordinates {
    (Laptop, 65.7528)
    (Workstation 1, 75.9762)
    (Workstation 2, 27.9005)
    (Server, 77.5292)
};

\addplot+[pattern color=nvidiaGreen, pattern=north east lines, draw=nvidiaGreen] coordinates {
    (Laptop, 6.91417)
    (Workstation 1, 11.7452)
    (Workstation 2, 1.82841)
    (Server, 3.34757)
};

\legend{CPU, GPU}

\end{axis}
\end{tikzpicture}
    \caption{Poisson on the cube ($p=5$). Matrix size: 140{,}556 for Laptop, Workstation~1, and Server; 21{,}924 for Workstation~2 (limited by VRAM).}
    \end{subfigure}

    \begin{subfigure}{0.9\columnwidth}
\begin{tikzpicture}
\begin{axis}[
    ybar,
    bar width=12pt,
    width=0.9\columnwidth,
    height=4cm,
    enlarge x limits={0.2},
    ylabel={Time to solve (s)},
    symbolic x coords={Laptop,Workstation 1,Server,Workstation 2},
    xtick=data,
    legend style={
        at={(1.02,0.5)},
        anchor=west
    },
    nodes near coords,
    nodes near coords align={vertical},
    ymajorgrids=true
]

\addplot+[fill=blue!50] coordinates {
    (Laptop, 51.319)
    (Workstation 1, 61.1513)
    (Workstation 2, 9.48666)
    (Server, 63.2592)
};

\addplot+[pattern color=nvidiaGreen, pattern=north east lines, draw=nvidiaGreen] coordinates {
    (Laptop, 5.90079)
    (Workstation 1, 9.44469)
    (Workstation 2, 1.34033)
    (Server, 2.64741)
};

\legend{CPU, GPU}

\end{axis}
\end{tikzpicture}
    \caption{Poisson on the L-shaped domain ($p=4$). Matrix size: 385{,}050 for Laptop, Workstation~1, and Server; 55{,}458 for Workstation~2 (limited by VRAM).}
    \end{subfigure}
    \begin{subfigure}{0.9\columnwidth}
\begin{tikzpicture}
\begin{axis}[
    ybar,
    bar width=12pt,
    width=0.9\columnwidth,
    height=4cm,
    enlarge x limits={0.2},
    ylabel={Time to solve (s)},
    symbolic x coords={Laptop,Workstation 1,Server,Workstation 2},
    xtick=data,
    legend style={
        at={(1.02,0.5)},
        anchor=west
    },
    nodes near coords,
    nodes near coords align={vertical},
    ymajorgrids=true
]

\addplot+[fill=blue!50] coordinates {
    (Laptop, 22.7514)
    (Workstation 1, 26.5754)
    (Workstation 2, 29.3357)
    (Server, 26.9735)
};

\addplot+[pattern color=nvidiaGreen, pattern=north east lines, draw=nvidiaGreen] coordinates {
    (Laptop, 7.48243)
    (Workstation 1, 4.13956)
    (Workstation 2, 4.29882)
    (Server, 1.23670)
};

\legend{CPU, GPU}

\end{axis}
\end{tikzpicture}
    \caption{Poisson on the quarter-ring domain ($p=4$). Matrix size: 132{,}600 for all systems.}
    \end{subfigure}
    
    \caption{Solve times on different systems (CPU vs GPU) for the three benchmark problems.}
    \label{fig:consumer_hardware}
\end{figure}
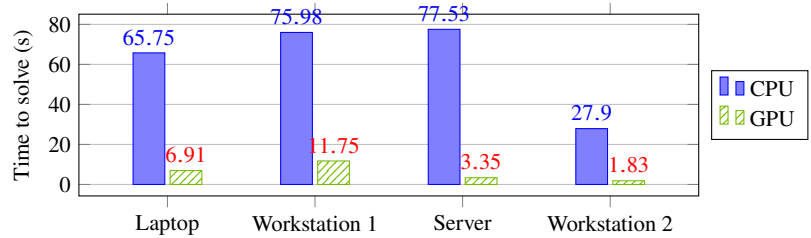
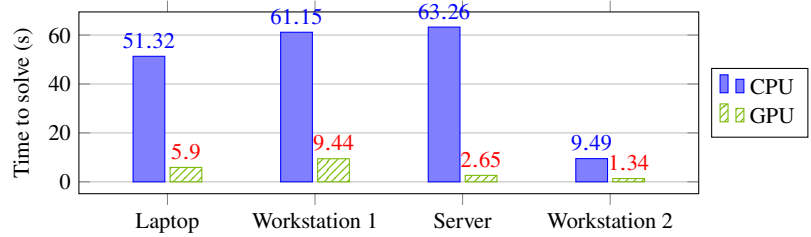
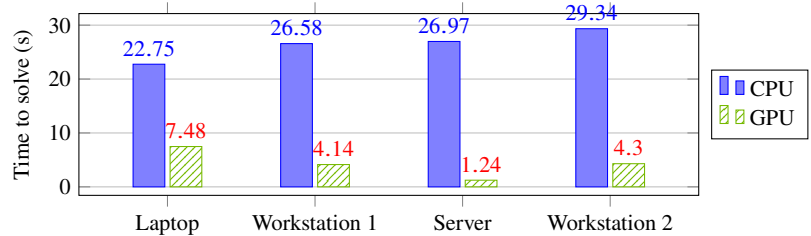

Across all platforms and test problems, the GPU implementation provides a clear and consistent reduction in solution time. For the Poisson cube case, speedups range from approximately $9.5\times$ on the Laptop to $23\times$ on the Server, with Workstation~1 and Workstation~2 achieving $6.5\times$ and $15\times$, respectively. For the L-shaped domain, the gains are similarly strong, ranging from $8.7\times$ (Laptop) and $6.5\times$ (Workstation~1) to $23.9\times$ (Server) and $7.1\times$ (Workstation~2). For the quarter-ring geometry, GPU acceleration yields speedups between $3.0\times$ (Laptop) and $22\times$ (Server), with intermediate improvements of $6.4\times$ and $6.8\times$ on the two workstations.

These results demonstrate that even mid-range consumer GPUs, such as the RTX~4060 in the laptop configuration, deliver order-of-magnitude speedups over their CPU counterparts. Higher-end professional GPUs, such as the NVIDIA~A40, further amplify this advantage, achieving more than $20\times$ acceleration on several test cases. Overall, these findings highlight the ability of \texttt{PSCToolkit} to deliver efficient GPU-based solvers across a broad range of hardware, extending high-performance numerical simulation capabilities well beyond traditional HPC infrastructures.

\section{Conclusions and future developments}
\label{sec:conclusions}

We have demonstrated that AMG-type preconditioners based on compatible weighted matching, as implemented in the \texttt{PSCToolkit} library suite~\cite{DAMBRA2023100463}, effectively accelerate the convergence of the conjugate gradient method for linear systems arising from the isogeometric discretization of elliptic equations. When combined with standard smoothers, the resulting solvers exhibit scalability with respect to both mesh refinement and the number of processing units, making them suitable for large-scale parallel simulations.

A promising direction for future work is the development of hybrid geometric--algebraic multigrid strategies tailored to problems of this type. In such an approach, a geometric multigrid scheme would first reduce the problem dimensionality, while the coarse-level correction would be carried out using the algebraic multigrid preconditioner introduced in this work. This combination has the potential to exploit geometric structure while retaining the robustness of algebraic techniques on coarse levels. Another avenue of research concerns extending these methods to locally refined IgA discretizations, such as those based on truncated hierarchical B-splines~\cite{MR3711991} or LR-splines~\cite{MR4118824}. In this setting, purely geometric multigrid approaches often struggle to define consistent restriction and prolongation operators across locally refined regions{, even though approaches working under the assumption that the hierarchical mesh satisfies certain admissibility conditions are already available in the literature~\cite{MR4454927}.} Hybrid geometric--algebraic strategies, however, could naturally accommodate such nonuniform refinement patterns by leveraging geometric information where available and employing algebraic constructions to preserve stability and efficiency across irregularly refined levels. Such methods would therefore provide a flexible and scalable framework for addressing complex geometries and adaptively refined domains within the IgA paradigm.

In addition, we plan to extend GPU support in \texttt{PSCToolkit} through the integration of \texttt{OpenACC}, a directive-based programming model that enables the acceleration of existing codes on heterogeneous architectures, most notably GPUs, through simple compiler annotations. Unlike low-level programming models such as CUDA, OpenACC allows developers to introduce parallelism incrementally while maintaining a unified code base that remains portable across CPUs and different GPU vendors, thereby extending support to AMD GPUs, which are increasingly widespread in TOP500 systems\footnote{See the sublist generator at \href{https://top500.org/statistics/list/}{top500.org/statistics/list/}. In the November~2025 ranking, NVIDIA hardware accounts for 43.8\% of the system share, compared to AMD's 5.8\%. In terms of performance share, AMD holds 37.8\%, while NVIDIA is at 51.3\%.}. Incorporating OpenACC will further enhance the portability and usability of \texttt{PSCToolkit}, enabling a broader range of users to benefit from GPU acceleration with minimal code modifications.

\begin{acknowledgement}
FD acknowledges the MUR Excellence Department Project awarded to the Department of Mathematics, University of Pisa, CUP I57G22000700001. All authors are members of INdAM-GNCS. This work was partially supported by: Spoke 1 ``FutureHPC \& BigData'' and Spoke 6 ``Multiscale Modelling \& Engineering Applications'' of the Italian Research Center on High-Performance Computing, Big Data and Quantum Computing (ICSC) funded by MUR Missione 4 Componente 2 Investimento 1.4: Potenziamento strutture di ricerca e creazione di ``campioni nazionali di R\&S (M4C2-19)'' - Next Generation EU (NGEU); and by the INdAM-GNCS Project CUP E53C24001950001.

We thank the two referees for their valuable comments and suggestions, which have greatly improved the presentation of this work.
\end{acknowledgement}
\ethics{Competing Interests}{The authors have no conflicts of interest to declare that are relevant to the content of this chapter.}

\bibliographystyle{spmpsci}
\bibliography{bibliga}

\begin{thebibliography}{10}
\providecommand{\url}[1]{{#1}}
\providecommand{\urlprefix}{URL }
\expandafter\ifx\csname urlstyle\endcsname\relax
  \providecommand{\doi}[1]{DOI~\discretionary{}{}{}#1}\else
  \providecommand{\doi}{DOI~\discretionary{}{}{}\begingroup
  \urlstyle{rm}\Url}\fi

\bibitem{Quad2}
Antolin, P., Buffa, A., Calabr\`o, F., Martinelli, M., Sangalli, G.: Efficient
  matrix computation for tensor-product isogeometric analysis: the use of sum
  factorization.
\newblock Comput. Methods Appl. Mech. Engrg. \textbf{285}, 817--828 (2015).
\newblock \doi{10.1016/j.cma.2014.12.013}.
\newblock \urlprefix\url{https://doi.org/10.1016/j.cma.2014.12.013}

\bibitem{Arnold20061}
Arnold, D.N., Falkt, R.S., Whither, R.: Finite element exterior calculus,
  homological techniques, and applications.
\newblock Acta Numerica \textbf{15}, 1 – 155 (2006).
\newblock \doi{10.1017/S0962492906210018}

\bibitem{MR2861652}
Baker, A.H., Falgout, R.D., Kolev, T.V., Yang, U.M.: Multigrid smoothers for
  ultraparallel computing.
\newblock SIAM J. Sci. Comput. \textbf{33}(5), 2864--2887 (2011).
\newblock \doi{10.1137/100798806}.
\newblock \urlprefix\url{https://doi.org/10.1137/100798806}

\bibitem{Sangalli2}
Bosy, M., Montardini, M., Sangalli, G., Tani, M.: A domain decomposition method
  for isogeometric multi-patch problems with inexact local solvers.
\newblock Comput. Math. Appl. \textbf{80}(11), 2604--2621 (2020).
\newblock \doi{10.1016/j.camwa.2020.08.024}.
\newblock \urlprefix\url{https://doi.org/10.1016/j.camwa.2020.08.024}

\bibitem{MR3951499}
Bressan, A., Takacs, S.: Sum factorization techniques in isogeometric analysis.
\newblock Comput. Methods Appl. Mech. Engrg. \textbf{352}, 437--460 (2019).
\newblock \doi{10.1016/j.cma.2019.04.031}.
\newblock \urlprefix\url{https://doi.org/10.1016/j.cma.2019.04.031}

\bibitem{MR2808112}
Buffa, A., de~Falco, C., Sangalli, G.: Iso{G}eometric {A}nalysis: stable
  elements for the 2{D} {S}tokes equation.
\newblock Internat. J. Numer. Methods Fluids \textbf{65}(11-12), 1407--1422
  (2011).
\newblock \doi{10.1002/fld.2337}.
\newblock \urlprefix\url{https://doi.org/10.1002/fld.2337}

\bibitem{BPXBuffa}
Buffa, A., Harbrecht, H., Kunoth, A., Sangalli, G.: B{PX}-preconditioning for
  isogeometric analysis.
\newblock Comput. Methods Appl. Mech. Engrg. \textbf{265}, 63--70 (2013).
\newblock \doi{10.1016/j.cma.2013.05.014}.
\newblock \urlprefix\url{https://doi.org/10.1016/j.cma.2013.05.014}

\bibitem{Quad3}
Calabr\`o, F., Loli, G., Sangalli, G., Tani, M.: Quadrature rules in the
  isogeometric {G}alerkin method: state of the art and an introduction to
  weighted quadrature.
\newblock In: Advanced methods for geometric modeling and numerical simulation,
  \emph{Springer INdAM Ser.}, vol.~35, pp. 43--55. Springer, Cham (2019)

\bibitem{6009071}
{\c{C}}ataly{\"u}rek, {\"U}.V., Dobrian, F., Gebremedhin, A., Halappanavar, M.,
  Pothen, A.: Distributed-memory parallel algorithms for matching and coloring.
\newblock In: 2011 IEEE International Symposium on Parallel and Distributed
  Processing Workshops and Phd Forum, pp. 1971--1980 (2011).
\newblock \doi{10.1109/IPDPS.2011.360}

\bibitem{MR4331965}
D'Ambra, P., Durastante, F., Filippone, S.: A{MG} preconditioners for linear
  solvers towards extreme scale.
\newblock SIAM J. Sci. Comput. \textbf{43}(5), S679--S703 (2021).
\newblock \doi{10.1137/20M134914X}.
\newblock \urlprefix\url{https://doi.org/10.1137/20M134914X}

\bibitem{DAmbra2025}
D'Ambra, P., Durastante, F., Filippone, S., Massei, S., Thomas, S.: Optimal
  polynomial smoothers for parallel {AMG}.
\newblock Numerical Algorithms  (2025).
\newblock \doi{10.1007/s11075-025-02117-6}.
\newblock \urlprefix\url{https://doi.org/10.1007/s11075-025-02117-6}

\bibitem{DAmbra2018BootCMatch}
D'Ambra, P., Filippone, S., Vassilevski, P.S.: Boot{CM}atch: {A} {S}oftware
  {P}ackage for {B}ootstrap {AMG} {B}ased on {G}raph {W}eighted {M}atching.
\newblock ACM Transactions on Mathematical Software (TOMS) \textbf{44}(2),
  16:1--16:26 (2018).
\newblock \doi{10.1145/3190647}.
\newblock \urlprefix\url{https://doi.org/10.1145/3190647}

\bibitem{GeoPDEs1}
{de Falco}, C., Reali, A., Vázquez, R.: {GeoPDEs}: {A} research tool for
  {I}sogeometric {A}nalysis of {PDE}s.
\newblock Advances in Engineering Software \textbf{42}(12), 1020--1034 (2011).
\newblock \doi{https://doi.org/10.1016/j.advengsoft.2011.06.010}.
\newblock
  \urlprefix\url{https://www.sciencedirect.com/science/article/pii/S0965997811001839}

\bibitem{MR3644391}
Diestel, R.: Graph theory, \emph{Graduate Texts in Mathematics}, vol. 173,
  fifth edn.
\newblock Springer, Berlin (2017).
\newblock \doi{10.1007/978-3-662-53622-3}.
\newblock \urlprefix\url{https://doi.org/10.1007/978-3-662-53622-3}

\bibitem{Manni1}
Donatelli, M., Garoni, C., Manni, C., Serra-Capizzano, S., Speleers, H.: Robust
  and optimal multi-iterative techniques for {I}g{A} {G}alerkin linear systems.
\newblock Comput. Methods Appl. Mech. Engrg. \textbf{284}, 230--264 (2015).
\newblock \doi{10.1016/j.cma.2014.06.001}.
\newblock \urlprefix\url{https://doi.org/10.1016/j.cma.2014.06.001}

\bibitem{DAMBRA2023100463}
D’Ambra, P., Durastante, F., Filippone, S.: Parallel sparse computation
  toolkit.
\newblock Software Impacts \textbf{15}, 100463 (2023).
\newblock \doi{https://doi.org/10.1016/j.simpa.2022.100463}.
\newblock
  \urlprefix\url{https://www.sciencedirect.com/science/article/pii/S2665963822001476}

\bibitem{Evans2020422}
Evans, J.A., Scott, M.A., Shepherd, K.M., Thomas, D.C., Vázquez~Hernández,
  R.: Hierarchical b-spline complexes of discrete differential forms.
\newblock IMA Journal of Numerical Analysis \textbf{40}(1), 422 – 473 (2020).
\newblock \doi{10.1093/imanum/dry077}

\bibitem{Evans}
Evans, L.C.: Partial differential equations, \emph{Graduate Studies in
  Mathematics}, vol.~19.
\newblock American Mathematical Society, Providence, RI (1998).
\newblock \doi{10.1090/gsm/019}.
\newblock \urlprefix\url{https://doi.org/10.1090/gsm/019}

\bibitem{10.1145/3017994}
Filippone, S., Cardellini, V., Barbieri, D., Fanfarillo, A.: Sparse
  matrix-vector multiplication on gpgpus.
\newblock ACM Trans. Math. Softw. \textbf{43}(4) (2017).
\newblock \doi{10.1145/3017994}.
\newblock \urlprefix\url{https://doi.org/10.1145/3017994}

\bibitem{GenericAMG}
Gahalaut, K.P.S., Tomar, S.K., Kraus, J.K.: Algebraic multilevel
  preconditioning in isogeometric analysis: construction and numerical studies.
\newblock Comput. Methods Appl. Mech. Engrg. \textbf{266}, 40--56 (2013).
\newblock \doi{10.1016/j.cma.2013.07.002}.
\newblock \urlprefix\url{https://doi.org/10.1016/j.cma.2013.07.002}

\bibitem{MR3711991}
Garau, E.M., V\'azquez, R.: Algorithms for the implementation of adaptive
  isogeometric methods using hierarchical {B}-splines.
\newblock Appl. Numer. Math. \textbf{123}, 58--87 (2018).
\newblock \doi{10.1016/j.apnum.2017.08.006}.
\newblock \urlprefix\url{https://doi.org/10.1016/j.apnum.2017.08.006}

\bibitem{Hiemstra20141444}
Hiemstra, R., Toshniwal, D., Huijsmans, R., Gerritsma, M.: High order geometric
  methods with exact conservation properties.
\newblock Journal of Computational Physics \textbf{257}(PB), 1444 – 1471
  (2014).
\newblock \doi{10.1016/j.jcp.2013.09.027}

\bibitem{MR4454927}
Hofreither, C., Mitter, L., Speleers, H.: Local multigrid solvers for adaptive
  isogeometric analysis in hierarchical spline spaces.
\newblock IMA J. Numer. Anal. \textbf{42}(3), 2429--2458 (2022).
\newblock \doi{10.1093/imanum/drab041}.
\newblock \urlprefix\url{https://doi.org/10.1093/imanum/drab041}

\bibitem{MR3686804}
Hofreither, C., Takacs, S.: Robust multigrid for isogeometric analysis based on
  stable splittings of spline spaces.
\newblock SIAM J. Numer. Anal. \textbf{55}(4), 2004--2024 (2017).
\newblock \doi{10.1137/16M1085425}.
\newblock \urlprefix\url{https://doi.org/10.1137/16M1085425}

\bibitem{VuikStokes}
Horn{\'i}kov{\'a}, H., Vuik, C.: {P}reconditioning for {L}inear {S}ystems
  {A}rising from {IgA} {D}iscretized {I}ncompressible {N}avier--{S}tokes
  {E}quations.
\newblock In: H.~van Brummelen, C.~Vuik, M.~M{\"o}ller, C.~Verhoosel,
  B.~Simeon, B.~J{\"u}ttler (eds.) Isogeometric Analysis and Applications 2018,
  pp. 77--97. Springer International Publishing, Cham (2021)

\bibitem{HUGHES2021467}
Hughes, T.J., Sangalli, G., Takacs, T., Toshniwal, D.: Chapter 8 - smooth
  multi-patch discretizations in isogeometric analysis.
\newblock In: A.~Bonito, R.H. Nochetto (eds.) Geometric Partial Differential
  Equations - Part II, \emph{Handbook of Numerical Analysis}, vol.~22, pp.
  467--543. Elsevier (2021).
\newblock \doi{https://doi.org/10.1016/bs.hna.2020.09.002}.
\newblock
  \urlprefix\url{https://www.sciencedirect.com/science/article/pii/S1570865920300144}

\bibitem{Quad1}
Hughes, T.J.R., Reali, A., Sangalli, G.: Efficient quadrature for {NURBS}-based
  isogeometric analysis.
\newblock Comput. Methods Appl. Mech. Engrg. \textbf{199}(5-8), 301--313
  (2010).
\newblock \doi{10.1016/j.cma.2008.12.004}.
\newblock \urlprefix\url{https://doi.org/10.1016/j.cma.2008.12.004}

\bibitem{METIS}
Karypis, G., Kumar, V.: A fast and high quality multilevel scheme for
  partitioning irregular graphs.
\newblock SIAM J. Sci. Comput. \textbf{20}(1), 359--392 (1998).
\newblock \doi{10.1137/S1064827595287997}.
\newblock \urlprefix\url{https://doi.org/10.1137/S1064827595287997}

\bibitem{IETIfirst}
Kleiss, S.K., Pechstein, C., J\"uttler, B., Tomar, S.: I{ETI}---isogeometric
  tearing and interconnecting.
\newblock Comput. Methods Appl. Mech. Engrg. \textbf{247/248}, 201--215 (2012).
\newblock \doi{10.1016/j.cma.2012.08.007}.
\newblock \urlprefix\url{https://doi.org/10.1016/j.cma.2012.08.007}

\bibitem{Lottes}
Lottes, J.: Optimal polynomial smoothers for multigrid {V}-cycles.
\newblock Numer. Linear Algebra Appl. \textbf{30}(6), Paper No. e2518, 27
  (2023).
\newblock \doi{10.1002/nla.2518}.
\newblock \urlprefix\url{https://doi.org/10.1002/nla.2518}

\bibitem{Mazza1}
Mazza, M., Manni, C., Ratnani, A., Serra-Capizzano, S., Speleers, H.:
  Isogeometric analysis for 2{D} and 3{D} curl-div problems: spectral symbols
  and fast iterative solvers.
\newblock Comput. Methods Appl. Mech. Engrg. \textbf{344}, 970--997 (2019).
\newblock \doi{10.1016/j.cma.2018.10.008}.
\newblock \urlprefix\url{https://doi.org/10.1016/j.cma.2018.10.008}

\bibitem{Sangalli1}
Montardini, M., Sangalli, G., Schneckenleitner, R., Takacs, S., Tani, M.: A
  {IETI}-{DP} method for discontinuous {G}alerkin discretizations in
  isogeometric analysis with inexact local solvers.
\newblock Math. Models Methods Appl. Sci. \textbf{33}(10), 2085--2111 (2023).
\newblock \doi{10.1142/S0218202523500495}.
\newblock \urlprefix\url{https://doi.org/10.1142/S0218202523500495}

\bibitem{Sangalli3}
Montardini, M., Sangalli, G., Tani, M.: Robust isogeometric preconditioners for
  the {S}tokes system based on the fast diagonalization method.
\newblock Comput. Methods Appl. Mech. Engrg. \textbf{338}, 162--185 (2018).
\newblock \doi{10.1016/j.cma.2018.04.017}.
\newblock \urlprefix\url{https://doi.org/10.1016/j.cma.2018.04.017}

\bibitem{FCG}
Notay, Y.: Flexible conjugate gradients.
\newblock SIAM J. Sci. Comput. \textbf{22}(4), 1444--1460 (2000).
\newblock \doi{10.1137/S1064827599362314}.
\newblock \urlprefix\url{https://doi.org/10.1137/S1064827599362314}

\bibitem{MR4087175}
Pan, M., J\"uttler, B., Giust, A.: Fast formation of isogeometric {G}alerkin
  matrices via integration by interpolation and look-up.
\newblock Comput. Methods Appl. Mech. Engrg. \textbf{366}, 113005, 25 (2020).
\newblock \doi{10.1016/j.cma.2020.113005}.
\newblock \urlprefix\url{https://doi.org/10.1016/j.cma.2020.113005}

\bibitem{MR4118824}
Patrizi, F., Manni, C., Pelosi, F., Speleers, H.: Adaptive refinement with
  locally linearly independent {LR} {B}-splines: theory and applications.
\newblock Comput. Methods Appl. Mech. Engrg. \textbf{369}, 113230, 20 (2020).
\newblock \doi{10.1016/j.cma.2020.113230}.
\newblock \urlprefix\url{https://doi.org/10.1016/j.cma.2020.113230}

\bibitem{zbMATH00792231}
Piegl, L., Tiller, W.: The {{\(NURBS\)}} book.
\newblock Berlin: Springer-Verlag (1995)

\bibitem{MR1990645}
Saad, Y.: Iterative methods for sparse linear systems, second edn.
\newblock Society for Industrial and Applied Mathematics, Philadelphia, PA
  (2003).
\newblock \doi{10.1137/1.9780898718003}.
\newblock \urlprefix\url{https://doi.org/10.1137/1.9780898718003}

\bibitem{Sangalli4}
Sangalli, G., Tani, M.: Isogeometric preconditioners based on fast solvers for
  the {S}ylvester equation.
\newblock SIAM J. Sci. Comput. \textbf{38}(6), A3644--A3671 (2016).
\newblock \doi{10.1137/16M1062788}.
\newblock \urlprefix\url{https://doi.org/10.1137/16M1062788}

\bibitem{Schumaker}
Schumaker, L.L.: Spline functions: basic theory, third edn.
\newblock Cambridge Mathematical Library. Cambridge University Press, Cambridge
  (2007).
\newblock \doi{10.1017/CBO9780511618994}.
\newblock \urlprefix\url{https://doi.org/10.1017/CBO9780511618994}

\bibitem{MR1393006}
Van\v{e}k, P., Mandel, J., Brezina, M.: Algebraic multigrid by smoothed
  aggregation for second and fourth order elliptic problems.
\newblock Computing \textbf{56}(3), 179--196 (1996).
\newblock \doi{10.1007/BF02238511}.
\newblock \urlprefix\url{https://doi.org/10.1007/BF02238511}.
\newblock International GAMM-Workshop on Multi-level Methods (Meisdorf, 1994)

\bibitem{GeoPDEs2}
V\'azquez, R.: A new design for the implementation of isogeometric analysis in
  {O}ctave and {M}atlab: {G}eo{PDE}s 3.0.
\newblock Comput. Math. Appl. \textbf{72}(3), 523--554 (2016).
\newblock \doi{10.1016/j.camwa.2016.05.010}.
\newblock \urlprefix\url{https://doi.org/10.1016/j.camwa.2016.05.010}

\bibitem{Veiga3}
Beir\~ao~da Veiga, L., Cho, D., Pavarino, L.F., Scacchi, S.: B{DDC}
  preconditioners for isogeometric analysis.
\newblock Math. Models Methods Appl. Sci. \textbf{23}(6), 1099--1142 (2013).
\newblock \doi{10.1142/S0218202513500048}.
\newblock \urlprefix\url{https://doi.org/10.1142/S0218202513500048}

\bibitem{Veiga2}
Beir\~ao~da Veiga, L., Cho, D., Pavarino, L.F., Scacchi, S.: Isogeometric
  {S}chwarz preconditioners for linear elasticity systems.
\newblock Comput. Methods Appl. Mech. Engrg. \textbf{253}, 439--454 (2013).
\newblock \doi{10.1016/j.cma.2012.10.011}.
\newblock \urlprefix\url{https://doi.org/10.1016/j.cma.2012.10.011}

\bibitem{Veiga}
Beir\~ao~da Veiga, L., Cho, D., Pavarino, L.F., Scacchi, S.: Overlapping
  {S}chwarz preconditioners for isogeometric collocation methods.
\newblock Comput. Methods Appl. Mech. Engrg. \textbf{278}, 239--253 (2014).
\newblock \doi{10.1016/j.cma.2014.05.007}.
\newblock \urlprefix\url{https://doi.org/10.1016/j.cma.2014.05.007}

\bibitem{XuZikatanov}
Xu, J., Zikatanov, L.: Algebraic multigrid methods.
\newblock Acta Numer. \textbf{26}, 591--721 (2017).
\newblock \doi{10.1017/S0962492917000083}.
\newblock \urlprefix\url{https://doi.org/10.1017/S0962492917000083}

\end{thebibliography}

\end{document}